\documentclass[11pt]{amsart}
\usepackage[english]{babel}
\usepackage{egothic}
\usepackage[T1]{fontenc}

\usepackage{amsmath}
\usepackage{amsthm}
\usepackage{amssymb}
\usepackage[all]{xy}
\usepackage{mathrsfs}
\usepackage{enumerate}
\usepackage{physics}
\usepackage[notcite, final, notref]{showkeys}
\usepackage{dsfont}
\usepackage[dvips]{color}
\newcommand{\scal}[1]{\left<#1\right>}
\newcommand{\Hq}{\mathbb H}
\newcommand{\Sq}{\mathbb S}

\newcommand{\N}{\mathbb{N}}

\newcommand{\R}{\mathbb{R}}      
\newcommand{\Z}{\mathbb{Z}}
\newcommand{\C}{\mathbb{C}}

\newcommand{\CFH}{\mathcal{F}_{CS}^\alpha(\mathbb{H})}

\newtheorem{theorem}{Theorem}[section]

\newtheorem{lemma}[theorem]{Lemma}
\newtheorem{proposition}[theorem]{Proposition}

\newtheorem{definition}[theorem]{Definition}

\theoremstyle{definition}
\newtheorem{remark}[theorem]{Remark}

\newtheorem{example}[theorem]{Example}

\renewcommand{\Re}{\mathrm{Re}}

\usepackage{xcolor}

\title[]{Towards Fock Spaces in Hypercomplex Analysis}
\author[K. Diki]{Kamal Diki}
\address{(KD) Clifford research group, Department of Electronics and Information Systems, Faculty of Engineering and Architecture, Ghent University, Krijgslaan 281, 9000 Ghent, Belgium}
\email{Kamal.Diki@UGent.be}

\begin{document}
\maketitle
\begin{abstract}
The Bargmann-Fock space(or Fock space for short) is a fundamental example of reproducing kernel Hilbert spaces that has found fascinating applications across multiple fields of current interest, including quantum mechanics, time-frequency analysis, mathematical analysis, and stochastic processes. In recent years, there has been increased interest in studying counterparts of the Fock space and related topics in hypercomplex analysis. This chapter presents a survey exploring various aspects of the Fock space from complex to hypercomplex analysis. In particular, we discuss different Fock spaces recently introduced in the setting of slice hyperholomorphic and slice polyanalytic functions of a quaternionic variable. The connection between slice hyperholomorphic (polyanalytic) Fock spaces and the classical theory of Fueter regular and poly-Fueter regular functions is established via the Fueter mapping theorem and
its polyanalytic extension. This chapter focuses on Fock spaces consisting of functions of a quaternionic
variable, with a brief discussion of related works in the Clifford setting.
\end{abstract}



\section{Introduction}
\setcounter{equation}{0}

In 1961, Bargmann introduced and studied an important Hilbert space of analytic functions and an associated integral transform in \cite{Barg}. This space is the so-called Segal-Bargmann-Fock space, also known in the literature as the bosonic Fock space \cite{Ner}. To simplify the presentation of this chapter, we refer to it simply as the Fock space. It consists of entire functions that are square integrable on the complex plane with respect to the normalized Gaussian measure.  The Fock space provides a natural setting for the creation and annihilation operators $$\displaystyle M_zf(z):=zf(z)\quad \text{and} \quad D_zf(z):=\frac{d}{dz}f(z),$$ which are closed, densely defined operators that are adjoints of each other in this space.  It is important to note that these operators satisfy the classical Heisenberg commutation rule $$\left[D_z,M_z\right]=\mathcal{I},$$ where $\left[\cdot,\cdot\right]$ denotes the commutator and $\mathcal{I}$ is the identity operator. The Fock space serves as a fundamental mathematical object in quantum mechanics, as it is unitarily equivalent to the classical $L^2$-Schrödinger Hilbert space on the real line via the Segal-Bargmann transform. This integral transform has found various applications across different areas of physics, including path integrals, coherent states, quantum field theory, and more recently, quantum gravity \cite{Gaz, Hall2013,SAD2024}. \\ \\

Fock spaces in the Banach case have been studied in mathematical analysis from different perspectives of complex analysis and operator theory. For example, see the classical book by Zhu \cite{Zhu2012} and the references therein. Moreover, various extensions of the Fock space in one complex variable in the Hilbert case have been explored. In \cite{Cholewinski1984} (see also \cite{SifiSoltani2002}), a generalized Fock space of even entire functions was investigated by choosing a weighted Hilbert space involving a modified Bessel function instead of the Gaussian measure. Another interesting extension of the Fock space in one complex variable is built upon the theory of polyanalytic functions. These functions were introduced in 1908 by Kolossov to solve problems in elasticity theory \cite{kolossov}. The theory of polyanalytic functions extends classical holomorphic functions by considering null-solutions of powers of the classical Cauchy-Riemann operator. For a general introduction to this topic, we refer to \cite{Balk1991, balk_ency}. In more recent times, this function theory has been studied by several authors from different perspectives \cite{Asampling, AbreuFeichtinger2014, AG2012, ACDS2022,ACDSS2023, Vas}. Note that the Fock space of polyanalytic functions can be introduced as the space of all polyanlytic functions on the complex plane that are square integrable with respect to the Gaussian measure \cite{AbreuFeichtinger2014,Balk1991}. \\ \\
 In hypercomplex analysis, two function theories are prominent: the recent theory of slice hyperholomorphic functions introduced in \cite{GentiliStruppa07}, and the classical theory of Fueter regular (monogenic) functions\cite{BDS1982,CSSS2004, GHS}. Slice hyperholomorphic functions were extended to the Clifford setting by considering slice monogenic functions in \cite{CSS2009}. This theory has been extensively developed in recent years and has found applications in Schur analysis \cite{SCHURBOOK}, quaternionic operator theory \cite{ColomboSabadiniStruppa2011}, and approximation theory \cite{GS2019}. It also provides a formalism for quaternionic quantum mechanics \cite{MTS2017}. This theory has been recently extended to the polyanalytic setting in \cite{ACDS20222, Rsimmath1, Rsimmath2, ADSP2020}. \\ \\  The theory of quaternionic Fueter regular (monogenic in the Clifford case) functions is defined by means of an extension of the Cauchy-Riemann operator in $ \mathbb{R}^4$, leading to the so-called Cauchy-Fueter operator (Dirac operator) \cite{BDS1982,CSSS2004,GHS}. A key difference between slice hyperholomorphic and Fueter regular functions is that polynomials and power series of the form $$\displaystyle P_N(q)=\sum_{n=0}^{N} q^na_n, \quad \text{ and} \quad f(q)=\sum_{n=0}^{\infty} q^na_n,$$ with quaternionic coefficients to the right (or to the left) are slice hyperholomorphic but not Fueter regular. The connection between the two function theories is established by the so-called Fueter mapping theorem, see \cite{ColomboSabadiniSommen2010,SCEBOOK}.  
 \\ \\The Fueter (or Fueter-Sce) mapping theorem is a fundamental result in hypercomplex analysis that allows the construction of Fueter regular functions starting from holomorphic functions. This result was originally proved in the quaternionic setting by Fueter in \cite{Fueter1}. The construction of the Fueter mapping theorem proceeds in two steps. First, slice hyperholomorphic functions are obtained by extending complex holomorphic functions using the slice extension. Then, the Laplace operator is applied to the slice hyperholomorphic extension to obtain Fueter regular functions. In the late 1950s, Sce extended the quaternionic Fueter mapping theorem to Clifford analysis in the case of odd dimensions \cite{SCEBOOK, S}. Later, Qian showed in \cite{Q} that Sce's theorem also holds in the case of even dimensions by employing Fourier multipliers. For further details and extensions of the Fueter mapping in hypercomplex analysis, we refer to \cite{ColomboSabadiniSommen2010, Qian2014,Q}. We also refer the reader to the recent book \cite{SCEBOOK}, which presents translations of Sce's works on this topic along with commentary.
   \\ \\
In recent years, the theory of Fock spaces has been extensively studied in hypercomplex analysis from various perspectives. The case of quaternionic variables, and more generally, the Clifford setting presents additional challenges compared to the classical complex case, mainly due to the lack of commutativity and the existence of multiple function theories extending the classical theory of holomorphic functions. For instance, the quaternionic slice hyperholomorphic Fock space and its associated Segal-Bargmann transform were introduced in \cite{ACSS2014,DG2017}. These constructions were applied to introduce and study a quaternionic counterpart of the Short-time Fourier transform(STFT) in \cite{DMD}. Moreover, an uncertaintly principle in the quaternionic Fock space of slice hyperholomorphic functions was established in \cite{Ren}.
Additionally, quaternionic (Banach) Fock spaces of the first and second kinds were introduced in \cite{DGS2019} to study approximation results. Subsequently, weighted compostion operators on these quaternionic Fock spaces were defined and studied in \cite{LL2025, LL2024, LL2021, LLL2023}. A quaternionic counterpart of the Mittag-Leffler Fock space and a Fueter-type Mittag-Leffler Bargmann transform were recently investigated in \cite{ADD2025}. On the other hand, Fock spaces of slice polyanalytic functions of a quaternionic variable have been considered in \cite{Rsimmath1, Rsimmath2, BEG, BBG2025, GS2024} for different purposes. Furthermore, a slice polyanalytic counterpart of the Short-time Fourier transform was developed in \cite{DMD2}, and generalized Appell polynomials and Fueter–Bargmann transforms in the polyanalytic setting were investigated in \cite{DMD2022}.\\ \\

\textit{Outline:} In Section 2, we recall basic definitions and properties of the Fock space in one complex variable, both in the analytic and polyanalytic cases. In Section 3, we introduce the standard notations of quaternions and review essential definitions of the quaternionic function theories relevant to this chapter. Specifically, we cover slice hyperholomorphic, slice polyanalytic, Fueter regular, and poly-Fueter regular function theories. Section 4 is devoted to studying several quaternionic Fock spaces of slice hyperholomorphic functions that have been considered in recent years. We also discuss some of their key properties, such as the reproducing kernel property and approximation results. In Section 5, we present quaternionic Fock spaces of slice polyanalytic functions. In Section 6, we recall the classical Fueter mapping theorem and introduce its polyanatic extensions. These theorems are used to construct Fock type spaces of Fueter and poly-Fueter regular functions starting from the slice hyperholomorphic and slice polyanalytic counterparts. Finally, in Section 7 we conclude the chapter by summarizing the main points discussed and suggesting new research avenues for Fock spaces in hypercomplex analysis.
\section{Fock spaces in one complex variable}
\setcounter{equation}{0}
In this section, we begin by reviewing basic definitions and properties of the Fock space. For further explanations, we refer to \cite{Folland,Hall2013,Ner,Zhu2012} in the analytic case and to \cite{AbreuFeichtinger2014, Balk1991, Vas} in the polyanalytic case.
\subsection{Analytic Case}

First, we present the geometric description of the Fock space, viewed as the subspace of entire functions on the complex plane that are square-integrable with respect to the Gaussian measure:

\begin{definition}\label{Fockdef}
An entire function $f: \mathbb{C} \to \mathbb{C}$ belongs to the Fock space, denoted by $ \mathcal{F}(\mathbb{C})$, if
$$ \| f \|_{\mathcal{F}(\mathbb{C})}^2= \frac{1}{\pi} \int_{\mathbb{C}} | f(z)|^2 e^{-|z|^2} d \lambda(z) < \infty,$$
where $ d \lambda(z)=dx dy$ is the classical Lebesgue measure on $\mathbb{C}$ for $z=x+iy$.
\end{definition}
The space $\mathcal{F}(\mathbb{C})$ is equipped with the inner product

\begin{equation}
\langle f,g \rangle_{\mathcal{F}(\mathbb{C})}=\frac{1}{\pi}\int_\mathbb{C}\overline{g(z)}f(z)e^{-|z|^2} d \lambda(z),
\end{equation}
for all $f,g \in  \mathcal{F}(\mathbb{C})$. 

\begin{remark}\label{ProbaMeasure}
Observe that the measure  $$d\mu(z):=\frac{1}{\pi}e^{-|z|^2} d \lambda(z),$$ defines a probability measure on $\mathbb{C}$. Indeed, we have

\begin{align*}
\displaystyle \int_\mathbb C d\mu(z)&=\frac{1}{\pi}\int_\mathbb{C} e^{-|z|^2} d \lambda(z)=1.\\
\end{align*}
\end{remark}



\begin{theorem} \label{SCFock}
The Fock space $\mathcal{F}(\mathbb{C})$ admits the following analytic characterization
\begin{equation}
	\label{seq }
	\mathcal{F}(\mathbb{C})= \left \{ f(z)=\sum_{k=0}^\infty z^k a_k \text{ }\mid  \text{ } (a_k)_{k \in \mathbb{N}_0} \subset \mathbb{C}, \quad \sum_{k=0}^\infty k!|a_k|^2  < \infty \right\}.
\end{equation}
Moreover, for $f(z)=\displaystyle \sum_{k=0}^{\infty} a_kz^k$ and $g(z)=\displaystyle  \sum_{k=0}^{\infty} b_kz^k $ in $\mathcal{F}(\mathbb{C})$, the inner product is given by

\begin{equation}
\displaystyle \langle f,g \rangle_{\mathcal{F}(\mathbb{C})}= \sum_{k=0}^{\infty} k!\overline{b_k} a_k. 
\end{equation}
\end{theorem}

The analytic characterization of the Fock space and the Cauchy-Schwarz inequality allow us to prove that functions belonging to the Fock space $\mathcal{F}(\mathbb{C})$ satisfy the following estimate.

\begin{lemma}\label{FockEstimate}
  Let $f\in \mathcal{F}(\mathbb{C})$. Then, the following inequality holds

\begin{equation}
|f(z)|\leq e^{\frac{|z|^2}{2}} \| f \|_{\mathcal{F}(\mathbb{C})}, \qquad \forall z\in\mathbb{C}.
\end{equation}
\end{lemma}

Using Lemma \ref{FockEstimate}, we see that all the evaluation maps are bounded on the Fock space. Therefore, by the Riesz representation theorem, one can prove that $\mathcal{F}(\mathbb{C})$ is a reproducing kernel Hilbert space. More precisely, we have the following result.
\begin{theorem}
The Fock space is a reproducing kernel Hilbert space, with the reproducing kernel given by
\begin{equation}
K(z,w)=e^{z \overline{w}}, \qquad \forall z,w \in \mathbb{C}.
\end{equation}
Moreover, the reproducing kernel property ensures that for every $f\in\mathcal{F}(\mathbb{C})$ and $w\in \mathbb{C} $, the value of $f(w)$ can be expressed as
\begin{equation}
f(w)=\displaystyle  \frac{1}{\pi}\int_{\mathbb{C}} \overline{K(z,w)} f(z) e^{-|z|^2}d \lambda(z).
\end{equation}

\end{theorem}

A key role in the Fock space theory is played by the functions $(e_k)_{k\in \mathbb{N}_0}$ obtained by normalizing the monomials. These functions are defined by
 \begin{equation}
  e_k(z):=\frac{z^k}{\sqrt{k!}}, \quad  z\in \mathbb{C}, \quad k\in\mathbb{N}_0.
 \end{equation}

  The normalized monomials $(e_k)_{k\in \mathbb{N}_0}$ form an essential building block for understanding the structure of the Fock space. For instance, the counterpart of the Zaremba formula for reproducing kernel Hilbert spaces is established in the following result.
\begin{theorem}\label{ZarembaFock}
The system $\lbrace e_k(z) \text{ }|\text{ } k\in \mathbb{N}_0\rbrace $ forms an orthonormal basis of the Fock space  $\mathcal{F}(\mathbb{C})$. Furthermore, the reproducing kernel is expressed using the Zaremba formula as
\begin{equation}
\displaystyle K(z,w)=\sum_{k=0}^{\infty}e_k(z)\overline{e_k(w)},
\end{equation}
for all $z,w\in \mathbb{C}$.
\end{theorem}

Consider now the Hermite polynomials, which can be introduced using their generating function

 \begin{equation}
\displaystyle  e^{2xz-z^2}=\sum_{k=0}^{\infty} \frac{H_k(x)}{k!}z^k,
\end{equation}
for every $z\in\mathbb{C}$ and $x\in \mathbb{R}$. In this representation, the Hermite polynomials appear as the Taylor coefficients of the analytic function $e^{2xz-z^2}$.

\begin{remark}
It is important to note that Hermite polynomials $\lbrace{H_k \text{ }|\text{ }k\in\mathbb N_0}\rbrace$ form an orthogonal basis of the weighted Hilbert space $L^2(\mathbb{R}, e^{-x^2}dx)$. Moreover,  we have \begin{equation}
\|H_k\|_{L^2(\mathbb{R}, e^{-x^2}dx)}=\sqrt{k!2^k \sqrt{\pi}}, \quad k\in\mathbb N_0.
\end{equation}
\end{remark}

The Hermite functions $h_k$ are obtained by multiplying the Hermite polynomials $H_k$ with the Gaussian function $e^{-\frac{x^2}{2}}$. Namely, we have 

\begin{equation}
h_k(x)=e^{-\frac{x^2}{2}}H_k(x), \quad k\in \mathbb{N}_0.
\end{equation}

The Hermite functions $(h_k)_{k\geq 0}$ form an orthogonal basis of the classical Hilbert space $L^2(\mathbb{R})$. 
\begin{definition}
An orthonormal basis of $L^2(\mathbb{R})$ is given by the normalized Hermite functions, defined as
\begin{equation}
\psi_k(x):=\frac{h_k(x)}{\|h_k\| _{L^2}}=\frac{h_k(x)}{\sqrt{k!2^k \sqrt{\pi}}},\quad  k\in \mathbb N_0.
\end{equation}
\end{definition}
Using the generating function of Hermite polynomials, one can introduce and compute the \textit{Segal-Bargmann kernel} as follows.

\begin{definition}
The Segal-Bargmann kernel is defined by
$$\displaystyle   A(z,x)=A_z(x) := \sum_{k=0}^{\infty}\psi_k(x)\frac{\overline{z}^k}{\sqrt{k!}}=\pi^{-1/4} e^{-\frac{1}{2}(x^2+\overline{z}^2)+\sqrt{2}\overline{z}x},$$
for every $z\in \mathbb{C}$ and $x\in\mathbb{R}$.
\end{definition}
\begin{remark}
The reproducing kernel of the Fock space can be factorized in terms of the inner product in $L^2(\mathbb{R})$ using the Segal-Bargmann kernel, leading to

$$\langle A_w, A_z \rangle_{L^2(\mathbb{R})}=e^{z\overline{w}}=K(z,w), \quad \forall z,w\in \mathbb{C}.$$
\end{remark}

The Segal-Bargmann transform, introduced in \cite{Barg}, is defined as follows:

\begin{definition}
The Segal-Bargmann transform of a function $\psi\in L^2(\mathbb{R})$ is given by

\begin{align*}
\displaystyle (\mathcal{B}\psi)(z)&:=\langle \psi, A_{z} \rangle_{L^2(\mathbb{R})} \\
&=\pi^{-1/4} \int_{\mathbb{R}} e^{-\frac{1}{2}(x^2+z^2)+\sqrt{2}zx}\psi(x)dx,
\end{align*}
where $z\in \mathbb{C}$ and $A_z$ is the Segal-Bargmann kernel.
\end{definition}

The Segal-Bargmann transform plays a central role in relating the Schrödinger and Bargmann-Fock representations of the Heisenberg group (see \cite{GK2001}). This connection is based on the following fundamental result established in \cite{Barg}.

\begin{theorem}
The Segal-Bargmann transform $\mathcal{B}$ defines a unitary operator that maps the Hilbert space $L^2(\mathbb{R})$ onto the Bargmann-Fock space $\mathcal{F}(\mathbb{C})$. 
\end{theorem}

More specifically, we make the following observation:
\begin{remark}

The integral transform $\mathcal{B}$ maps the normalized Hermite functions $\psi_k(x)$ onto the orthonormal basis of the Fock space $e_n(z)$. This action is given by
\begin{equation}
(\mathcal{B}\psi_k)(z)=\frac{z^k}{\sqrt{k!}}:=e_k(z), \quad k\in \mathbb N_0.
\end{equation}

\end{remark}

We conclude this discussion on analytic Fock spaces by recalling the $p$-Fock spaces, as described in \cite[p.36]{Zhu2012}.

\begin{definition} Let $0 <p < \infty$ and $\alpha>0$. The $p$-Fock space $F^{p}_{\alpha}(\mathbb{C})$ is defined as the space of all entire functions in $\mathbb{C}$ satisfying $$\displaystyle\frac{\alpha p}{2\pi}\int_{\mathbb{C}}\left |f(z)e^{-\alpha |z|^{2}|/2}\right |^{p}d \lambda(z)<+\infty,$$ where $d \lambda(z)=d x d y$, $z=x+i y$, is the area measure in the complex plane.
\end{definition}
\begin{remark}{\rm Endowed with 
$$\|f\|_{p, \alpha}=\displaystyle \left(\frac{\alpha p}{2\pi} \int_{\mathbb{C}}\left |f(z)e^{-\alpha |z|^{2}/2}\right |^{p}d \lambda(z) \right)^{1/p},$$
it is known that $F_{\alpha}^{p}$ is a Banach space for $1\le p <\infty$, and a complete metric space for
$\|\cdot\|^{p}_{p, \alpha}$ when $0<p<1$. Additionally, if $p=+\infty$, then  endowed with $$\|f\|_{\infty, \alpha}=\operatorname{ess}\sup\left\{|f(z)|e^{-\alpha |z|^{2}|/2}, \quad z\in \mathbb{C}\right\},$$ $F^{\infty}_{\alpha}$ is a Banach space.}
\end{remark}

\subsection{Polyanalytic Case}
We begin by introducing the notion of polyanalytic functions. For a general introduction to this function theory we refer to \cite{AbreuFeichtinger2014, Balk1991,balk_ency,VasBook} .
\begin{definition}
A complex valued function $f:\Omega\subset \mathbb{C}\longrightarrow \mathbb{C}$ of class $\mathcal{C}^n$ that belongs to the kernel of the $n$-th power ($n\geq 1$) of the classical Cauchy-Riemann operator $\displaystyle \frac{\partial}{\partial \overline{z}}$, meaning that $$\displaystyle \frac{\partial^n}{\partial \overline{z}^n}f(z)=0, \quad \forall  z\in\Omega,$$
 is called a polyanalytic function of order $n$. The space of all polyanalytic functions of order $n$ defined on a domain $\Omega$ is denoted by $H_n(\Omega)$. \end{definition}
\begin{remark}
 A complex valued function that is polyanalytic of order $n$ in $\mathbb{C}$ is called poly-entire of order $n$.
\end{remark}
An interesting fact about these functions is that any polyanalytic function of order $n$ can be decomposed in terms of $n$ analytic functions. Specifically, it can be expressed as
\begin{equation}\label{Polydeccomplex}
f(z)=\displaystyle \sum_{k=0}^{n-1}\overline{z}^kf_k(z),
\end{equation}
where each $f_k$ is an analytic function in $\Omega$. In particular, expanding each analytic component using the series expansion theorem leads to the representation

\begin{equation}\label{exp1}
f(z)=\displaystyle \sum_{k=0}^{n-1}\sum_{j=0}^{\infty}\overline{z}^kz^ja_{k,j},
\end{equation}
where $(a_{k,j})$ are complex coefficients.\\ \\

A Cauchy-type formula for polyanalytic functions appeared for the first
time in Th\'eodoresco's doctoral thesis \cite{teodoresco} and was later recalled
in \cite{balk_ency}:
\begin{theorem} [\cite{balk_ency}, Theorem 1.3]
If a function $f$ is  polyanalytic of order $n$ in a closed domain $\Omega$
bounded by
a rectifiable closed contour $\Gamma$, then the value of $f$ at any
point $z$ of the domain
$\Omega$ is expressed, using values of the function itself and its formal
derivatives at
points $t$ of the boundary $\Gamma$, by the formula
\begin{equation}\label{TEODOR}
f(z)=\frac{1}{2\pi i}\sum_{\ell=0}^{n-1}\int_{\Gamma} \frac{1}{\ell! (t-z)}
(\bar z-\bar t)^\ell \frac{\partial^\ell f}{\partial \bar t^\ell}\, dt.
\end{equation}
\end{theorem}
The formula contains a finite sum involving the $n$ kernels
$$
\pi_\ell(z,t)= \frac{1}{\ell! (t-z)} (\bar z-\bar t)^\ell,\ \ \ \ell=0,...,n-1.
$$
Another polyanalytic Cauchy formula, introduced in the quaternionic and Clifford settings (see \cite{brackx1976k, delanghe1978hypercomplex}), can be presented in the complex case as follows:
$$
f(z)=\int_{\partial \Omega} \sum_{\ell=0}^{n-1}(-2)^\ell L_\ell(w-z) \
\frac{\partial ^\ell}{\partial \bar{w}^\ell} f(w)\ dw,
$$
where $\partial \Omega$ is the boundary of the smooth bounded domain $\Omega$ in $\mathbb{C}$, $dw$
denotes the infinitesimal arc length, and 
$$
L_\ell(z):=\frac{1}{2\pi i \  z}\frac{({\rm Re}(z))^\ell}{\ell !}, \ \ \ \ell=0,...,n-1.
$$

\textbf{Polyanalytic Fock Space} \\

The polyanalytic Fock space of order $n$, discussed in \cite[pp.169-170]{Balk1991}, extends the classical analytic Fock space presented in Definition \ref{Fockdef} by considering poly-entire functions of order $n$ that are square-integrable with respect to the normalized Gaussian measure. Specifically, we define:
\begin{definition}
Let $n\in \mathbb N$. The polyanalytic Fock space of order $n$, denoted by $\mathcal{F}_n(\mathbb{C})$, is defined as
$$\mathcal{F}_n(\mathbb{C}):=\left\lbrace f\in H_n(\mathbb{C})\mid  \frac{1}{\pi}\int_{\mathbb{C}}|f(z)|^2e^{-|z|^2}d\lambda(z)<\infty \right\rbrace,$$
where $d\lambda(z)$ is the classical Lebesgue measure on the complex plane. 
\end{definition}
\begin{theorem}
The reproducing kernel associated with the polyanlaytic Fock space $\mathcal{F}_n(\mathbb{C})$ is given by
\begin{equation}\label{Kn}
K_n(z,w)=e^{z\overline{w}}\displaystyle \sum_{k=0}^{n-1}\frac{(-1)^k}{k!}{n \choose k+1}|z-w|^{2k},
\end{equation}
for every $z,w\in\mathbb{C}.$
\end{theorem}
\begin{remark}
 The polyanalytic Fock kernel \eqref{Kn} can be expressed in terms of the generalized Laguerre polynomials, leading to (see \cite{AbreuFeichtinger2014})
 
 \begin{equation}
K_n(z,w)=e^{z\overline{w}} L^{1}_{n-1}(|z-w|^2), \quad z, w \in \mathbb{C}. 
\end{equation}
Here, $L^{\beta}_{k}$ denotes the generalized Laguerre polynomials, defined as

\begin{equation}\label{Laguerre}
L^\beta_k(x):=\sum_{j=0}^{k}(-1)^j {k+\beta\choose k-j  } \frac{x^j}{j!}.
\end{equation}

\end{remark}
Finally, we recall the complex polyanalytic $p$-Fock spaces $\mathcal{F}^{p}_{\alpha, n}(\mathbb{C})$ introduced by Abreu and Gröchenig in \cite{AG2012}:

\begin{definition}\label{polyFp}
Let $0<p<\infty$, $\alpha>0$ and $n\in\mathbb{N}$. The polyanalytic $p$-Fock space $\mathcal{F}^{p}_{\alpha, n}(\mathbb{C})$ consists of all poly-entire functions of order $n$ such that 

$$\displaystyle ||f||_{p,\alpha}:=\left(\frac{\alpha p}{2\pi}\int_{\mathbb{C}}\left |f(z)e^{-\alpha |z|^{2}|/2}\right |^{p}\right)^{1/p}d \lambda(z)<+\infty.$$
\end{definition}

\section{Function theories of a quaternionic variable}
\setcounter{equation}{0}
The non-commutative algebra of quaternions $\mathbb H$ is defined as
$$\displaystyle \mathbb{H}:=\left\lbrace{q=x_0+x_1i+x_2j+x_3k \mid   \ x_0,x_1,x_2,x_3\in\mathbb{R}}\right\rbrace,$$ where the imaginary units $i, j,$ and $k$ satisfy the relations $$i^2=j^2=k^2=-1,$$ and 
$$ij=-ji=k,\quad  jk=-kj=i,\quad  ki=-ik=j.$$
The real and imaginary parts of a quaternion $q=x_0+x_1i+x_2j+x_3k$ are defined as 

$$\Re(q)=x_0,\quad \Im(q)= \vec{q}\,=x_1i+x_2j+x_3k,$$
while the conjugate and modulus of $q$ are given by
$$\overline{q}=\Re(q)-\Im(q),  \quad \vert{q}\vert=\sqrt{q\overline{q}}=\sqrt{x_0^2+x_1^2+x_2^2+x_3^2}.$$

Quaternionic conjugation satisfies the property $$\overline{ pq }= \overline{q}\, \overline{p},$$ for all $p,q\in \mathbb{H}$. \\ \\ The unit sphere of purely imaginary quaternions is denoted by $\mathbb S$ and defined as $$\mathbb S:= \displaystyle\left\lbrace{x_1i+x_2j+x_3k \mid x_1^2+x_2^2+x_3^2=1}\right\rbrace=\lbrace{q\in{\mathbb{H}} \mid q^2=-1}\rbrace.$$ Notice that any quaternion $q\in \mathbb{H}\setminus \mathbb{R}$ can be uniquely written as $q=x+I_q y$ for some real numbers $x$ and $y>0$, and some imaginary unit $I_q\in \mathbb{S}$. For a given $I\in{\mathbb{S}}$, one can define the slice $$\mathbb{C}_I:= \mathbb{R}+\mathbb{R}I.$$
Note that the slice $\mathbb{C}_I$ is isomorphic to the complex plane $\mathbb{C}$, so it can be considered as a complex plane in $\mathbb{H}$ passing through $0$, $1$ and $I$. The union of the slices $(\mathbb{C}_I)_{I\in \mathbb{S}}$ is the whole space of quaternions $$\displaystyle \mathbb{H}=\underset{I\in{\mathbb{S}}}{\bigcup}\mathbb{C}_I=\underset{I\in{\mathbb{S}}}{\bigcup}(\mathbb{R}+\mathbb{R}I).$$

An important class of domains on which slice hyperholomorphic functions are studied is the class of axially symmetric slice domains, defined as follows:
\begin{definition}
A domain $\Omega\subset \Hq$ is said to be a slice domain (or just an $s$-domain) if $\Omega\cap{\mathbb{R}}$ is nonempty and, for all $I\in{\mathbb{S}}$, the set $\Omega_I:=\Omega\cap{\C_I}$ is a domain of the complex plane $\C_I$.
Moreover, if for every $q=x+Iy\in{\Omega}$, the entire sphere $x+y\mathbb{S}:=\lbrace{x+Jy \mid \, J\in{\mathbb{S}}}\rbrace$
is contained in $\Omega$, we say that  $\Omega$ is an axially symmetric slice domain.
\end{definition}

\begin{definition}
Let $U\subseteq\mathbb{H}$ be an axially symmetric open set and $\mathcal{U} = \left\{ (x,y)\in\mathbb{R}^2\ \mid  x+ I y\subset U\right\}.$
A function $f:U\to \mathbb{H}$
 is called a left
 slice function, if it is of the form
 \[
 f(q) = \alpha(x,y) + I\beta(x,y)\qquad \text{for all } q = x + I y\in U,
 \]
with the two functions $\alpha, \beta: \mathcal{U}\to \mathbb{H}$ satisfying the compatibility conditions
$\alpha(x,-y) = \alpha(x,y)$, $\beta(x,-y) = -\beta(x,y)$.
\end{definition}


\subsection{Slice Hyperholomorphic and Slice Polyanalytic Functions}
In this section, we recall the basic definitions of slice hyperholomorphic and slice polyanalytic functions of a quaternionic variable. For further details on the theory of slice hyperholomorphic functions, we refer to \cite{SCHURBOOK, ColomboSabadiniStruppa2016, ColomboSabadiniStruppa2011, CSS2009, GS2019, GentiliStoppatoStruppa2013, GentiliStruppa07} and the references therein. Additional explanations on the theory of slice polyanalytic functions and related works can be found in \cite{ACDS20222, QHilbertSpaces, Rsimmath1, Rsimmath2, ADSP2020,BEG, DMD2}.

\begin{definition}
Let $\Omega$ be a domain in $\mathbb H$. A quaternionic valued function $f: \Omega \longrightarrow \mathbb{H}$ of class $\mathcal{C}^1$ is said to be (left) slice hyperholomorphic if, for every $I\in \mathbb{S}$, its restriction $f_I$ to the slice $\mathbb{C}_{I}$ is holomorphic on $\Omega_I := \Omega \cap \mathbb{C}$. Specifically, 
$$\overline{\partial_I} f(x+Iy):=
\dfrac{1}{2}\left(\frac{\partial }{\partial x}+I\frac{\partial }{\partial y}\right)f_I(x+yI)
$$
vanishes identically on $\Omega_I$ for every $I\in\mathbb{S}$. The set of slice hyperholomorphic functions will be denoted by $\mathcal{SH}(\Omega)$. 
\end{definition}

\begin{definition}
Let $\Omega$ be a domain in $\mathbb H$, and let $f$ be a slice hyperholomorphic function on $\Omega$.
The slice derivative $\partial_S f$ of $f$, is defined by
\begin{equation*}
\partial_S(f)(q):=
\left\{
\begin{array}{rl}
\partial_I(f)(q)& \text{if } q=x+yI, y\neq 0\\
\displaystyle\frac{\partial}{\partial{x}}(f)(x) & \text{if } q=x \text{ is real}.
\end{array}
\right.
\end{equation*}
\end{definition}

\begin{example}
The quaternionic polynomials $\displaystyle P_N(q)=\sum_{k=0}^{N}q^k a_k$, with quaternion coefficients on the right, are slice hyperholomorphic functions. 
\end{example}
The right quaternionic vector space of slice hyperholomorphic functions $\mathcal{SH}(\Omega)$ is endowed with the natural topology of uniform convergence on compact sets. Moreover, the characterization of such functions on a ball centered at the origin is given by the series expansion theorem proved in the original paper \cite{GentiliStruppa07}.

\begin{theorem}\label{expansion}
An $\mathbb{H}$-valued function $f$ is slice hyperholomorphic if and only if it has a series expansion of the form
$$f(q)=\sum_{n=0}^{+\infty} q^na_n$$
converging on $B(0,R):=\{q\in\mathbb{H}\mid |q|<R\}$.
\end{theorem}
An important fact about slice hyperholmorphic functions is that they satisfy the so-called \textit{Representation Formula}.

\begin{theorem}[Representation Formula]\label{repform}
Let $\Omega$ be an axially symmetric slice domain and $f\in{\mathcal{SH}(\Omega)}$. Then, for any $I,J\in{\mathbb{S}}$, we have the formula
$$
f(x+Jy)= \frac{1}{2}(1-JI)f_I(x+Iy)+ \frac{1}{2}(1+JI)f_I(x-Iy)
$$
for all $q=x+Jy\in{\Omega}$.
\end{theorem}

Another interesting approach, which allows slice hyperholomorphic functions to be defined as null-solutions of a certain global operator with non constant coefficients, was introduced in \cite{CGCS2013} (see also \cite{CS2014, P2019}). Namely, one can consider the following definition:
\begin{definition}
Let $\Omega$ be an open set in $\mathbb{H}$, and let $f:\Omega\longrightarrow \mathbb{H}$ a function of class  $\mathcal{C}^1$. We define the global operator $G(f)$ by
$$\displaystyle Gf(q):=|\vec{q}\,|^2\partial_{x_0}f(q)+\vec{q}\,\sum_{\ell=1}^3x_\ell\partial_{x_\ell}f(q),$$
for any $q=x_0+\vec{q} \in\Omega$.
\end{definition}
\begin{remark}
It was proved in \cite{CGCS2013} that any slice hyperholomorphic function belongs to the kernel of the global operator $G$ on axially symmetric slice domains.
\end{remark}

A modified version of the global operator $G$ was introduced in \cite{P2019} as $$\displaystyle Vf(q):=\partial_{x_0}f(q)+\frac{\vec{q}\,}{|\vec{q}\,|^2}\sum_{\ell=1}^3x_\ell \partial_{x_\ell}f(q), \quad  q\in\Omega\setminus \R.$$
\begin{remark}\label{GVR}
For suitable domains, the operators $G$ and $V$ are related by
$$Vf(q)=\frac{1}{|\vec{q}\,|^2}Gf(q), \qquad q\in \Omega\setminus \R.$$
If $V(f)$ admits a unique continuous extension to the entire domain $\Omega$, we shall denote this extension again by $V(f)$. For $n\in \mathbb N,$ the iterates $V^n(f)$ are defined inductively as follows: assume that $V^{n-1}(f)$ is of class $\mathcal{C}^1$ on $\Omega\setminus \mathbb{R}$, and that $V(V^{n-1}(f))$ admits a continuous extension to $\Omega$. We then set $$V^n(f):=V(V^{n-1}(f)),$$
where $V^n(f)$ again denotes this continuous extension on $\Omega.$
\end{remark}

Slice polyanalytic functions of a quaternionic variable can be introduced as an extension of slice hyperholomorphic functions by considering slice functions that belong to the kernel of the $n$-th power of the Cauchy-Riemann operator on every slice, as presented in the following definition.
\begin{definition}\label{Poly}
Let $\Omega$ be an axially symmetric open set in $\mathbb{H}$, and let $f: \Omega \longrightarrow \mathbb{H}$ be a left slice function $f(q)=f(x+Iy)=\alpha(x,y)+I \beta(x,y)$ of class  $\mathcal{C}^n(\Omega)$. Then, $f$ is called (left) slice polyanalytic of order $n$ on $\Omega$ if it satisfies on $\Omega$ the poly Cauchy-Riemann equations of order $n\in \mathbb{N}$, given by
$$
\overline{\partial_I}^n f(x+Iy):=
\dfrac{1}{2^n}\left(\frac{\partial }{\partial x}+I\frac{\partial }{\partial y}\right)^nf_I(x+yI)=0, \quad \text{ for all } I\in \mathbb{S}.
$$
\end{definition}
\begin{remark}
The set of all slice polyanalytic functions of order $n$ on an axially symmetric open set $\Omega$ is a right quaternionic vector space, denoted by $\mathcal{SP}_n(\Omega)$.
\end{remark}
The quaternionic counterpart of the complex polyanalytic decomposition \eqref{Polydeccomplex} for slice polyanalytic functions is given by the following result.

\begin{theorem}[Poly-decomposition]\label{carac3}
A function $f(q)$ is slice polyanalytic of order $n$ on a domain $\Omega$ if and only there exist  slice hyperholomorphic functions $f_0,...,f_{n-1}$ on $\Omega$ such that $$
f(q):=\displaystyle\sum_{k=0}^{n-1}\overline{q}^kf_k(q),
$$
 for every  $q\in\Omega$.
\end{theorem}

\begin{remark}
The Cauchy formulas corresponding to the theory of slice polyanalytic functions were introduced in \cite{ACDS20222,ADSP2020} to study polyanalytic counterparts of the the S-functional and the Fueter mapping theorem.
\end{remark}

\subsection{Fueter Regular and Poly-Fueter Regular Functions}

In this section, we recall the basic definitions of Fueter regular and poly-Fueter regular functions of a quaternionic variable. For further reading on the theory of Fueter regular (or monogenic) functions, we refer to \cite{BDS1982,CSSS2004, GHS} and the references therein. For the theory of poly-Fueter (or $k$-mongenic) functions, see \cite{brackx1976k, delanghe1978hypercomplex}. 

\begin{definition} Let $U\subset \mathbb{H}$ be an open set, and let $f:U\longrightarrow \mathbb{H}$ be a function of class $\mathcal{C}^1$. We say that $f$ is (left) Fueter regular on $U$ if it satisfies the equation $$D f(q):=\displaystyle\left(\frac{\partial}{\partial x_0}+i\frac{\partial}{\partial x_1}+j\frac{\partial}{\partial x_2}+k\frac{\partial}{\partial x_3}\right)f(q)=0, \text{ for all } q\in U,$$
where $D$ is known as the Cauchy-Fueter operator.
\end{definition}
\begin{remark}
The set of all Fueter regular functions on $U$ forms a right quaternionic vector space, denoted by $\mathcal{FR}(U)$.
\end{remark}
\begin{example}
The Fueter variables, defined by

$$\zeta_1(q)=x_1-x_0i, \quad \zeta_2(q)=x_2-x_0j, \quad \zeta_3(q)=x_3-x_0k, $$
 
 for all $q=x_0+x_1i+x_2j+x_3k\in \mathbb{H}$, are building block examples for the theory of Fueter regular functions.
\end{example}
 The poly-Fueter regular functions extend the notion of Fueter regular functions by considering functions that lied in the kernel of the $n$-th power of the Cauchy-Fueter operator $D$. These functions were introduced by Brackx in \cite{brackx1976k} for functions of a quaternionic variable and later extended to the Clifford setting by Brackx and Delanghe in \cite{delanghe1978hypercomplex}. 
\begin{definition}
Let $U\subset \mathbb{H}$ be an open set, and let $f:U\longrightarrow \mathbb{H}$ be a function of class $\mathcal{C}^n$. We say that $f$ is (left) poly-Fueter regular, or poly-regular for short, of order $n\geq 1$ on $U$ if it satisfies the equation $$D^n f(q):=\displaystyle\left(\frac{\partial}{\partial x_0}+i\frac{\partial}{\partial x_1}+j\frac{\partial}{\partial x_2}+k\frac{\partial}{\partial x_3}\right)^nf(q)=0, \text{ for all } q\in U.$$
\end{definition}

\begin{remark}
The set of all poly-Fueter regular functions of order $n$ on $U$ forms a right quaternionic vector space, denoted by $\mathcal{FR}_n(U)$.
\end{remark}

The proof of the following result appears in sections 6 and 7 of  \cite{brackx1976k} (see also \cite{delanghe1978hypercomplex} for the Clifford setting).
\begin{proposition}\label{Poly-Fueter-Dec}
A function $f$ is poly-Fueter regular of order $n$ if and only if it can be decomposed in terms of some unique Fueter regular functions $\phi_0,...,\phi_{n-1}$ such that $$f(q)=\displaystyle \sum_{k=0}^{n-1}x_0^k\phi_k(q),$$
 for every   $q=x_0+\vec{q} \in U$.
\end{proposition}

\section{Fock spaces of slice hyperholomorphic functions}
\setcounter{equation}{0}
In this section, we present different definitions of quaternionic Fock spaces and their corresponding reproducing kernels that have been studied in the slice hyperholomorphic setting.
\subsection{Quaternionic Slice Hyperholomorphic Fock space and Its Segal-Bargmann Transform}
The first quaternionic Fock space of slice hyperholomorphic functions has been introduced in \cite{ACSS2014} as follows:
\begin{definition}\label{SliceFock}
Let $I$ be any element in ${\mathbb{S}}$. The quaternionic (slice hyperholomorphic) Fock space is defined as
$$\displaystyle \mathcal{F}_{Slice}(\mathbb{H}):=\left\lbrace{f\in{\mathcal{SH}(\mathbb{H}) \mid \displaystyle \Vert{f}\Vert_{\mathcal{F}_{Slice}(\mathbb{H})}^2:=  \frac{1}{\pi} \int_{\C_I}\vert{f_I(p)}\vert^2 e^{-\vert{p}\vert^2}d\lambda_I(p) <\infty}}\right\rbrace,$$
 where $f_I = f|_{\C_I}$ is the restriction of $f$ to $\mathbb{C}_I$, and $d\lambda_I(p)=dxdy$ for $p=x+yI$.
 \end{definition}
The right quaternionic vector space $\mathcal{F}_{Slice}(\mathbb{H})$ is equipped with the inner product
\begin{equation}\label{spfg}
\langle f,g\rangle_{\mathcal{F}_{Slice}(\mathbb{H})} = \frac{1}{\pi}\int_{\C_I}\overline{g_I(p)}f_I(p)e^{-\vert{p}\vert^2} d\lambda_I(p),
\end{equation}
 where $f,g\in{\mathcal{F}_{Slice}(\mathbb{H})}$. 
     \begin{remark}
The definition of $\mathcal{F}_{Slice}(\mathbb{H})$ does not depend on the choice of the imaginary unit $I\in \mathbb{S}$ and this can be justified using the Representation Formula for slice hyperholomorphic functions presented in Theorem \ref{repform}. 
\end{remark}

The space $\mathcal{F}_{Slice}(\mathbb{H})$ contains the quaternionic monomials $\lbrace q^k \mid k\in \mathbb{N}_0\rbrace$ which form an orthogonal basis of this space. Indeed, we have 
 \begin{equation}\label{spmonomials}
 \langle q^k,q^s \rangle_{\mathcal{F}_{Slice}(\mathbb{H})}  = k!\delta_{k,s}.
\end{equation}

Additionally, it turns out that $\mathcal{F}_{Slice}(\mathbb{H})$ is a quaternionic reproducing kernel Hilbert space. Specifically, we have the following result:
\begin{theorem}
The slice hyperholomorphic Fock space $\mathcal{F}_{Slice}(\mathbb{H})$ is a right quaternionic reproducing kernel Hilbert space. Its reproducing kernel is given by the kernel function 
\begin{equation}\label{*Fockkernel}
K_q(p)=K(p,q):=e_*(p\bar{q})=\displaystyle \sum_{k=0}^\infty\frac{p^k\bar{q}^k}{k!},
\end{equation}
for every $(p,q)\in \mathbb{H}\times\mathbb{H}$. Moreover, for every $f\in \mathcal{F}_{Slice}(\mathbb{H})$ and $q\in \mathbb{H}$ we have 
$$\langle f, K_q  \rangle_{\mathcal{F}_{Slice}(\mathbb{H})}  =f(q).$$
\end{theorem}

A quaternionic counterpart of the Segal-Bargmann transform associated to the slice hyperholomorphic Fock space $\mathcal{F}_{Slice}(\mathbb{H})$ has been introduced and studied in \cite{DG2017}. This integral transform acts on the quaternionic Hilbert space $L^2_\mathbb{H}(\R)$ consisting of quaternionic valued functions $\psi:\R\longrightarrow \mathbb{H}$ that satisfy the $L^2$-integrability condition  $$ \displaystyle \Vert{\psi}\Vert_{L^2_\mathbb{H}(\R)}^2= \int_\R|\psi(x)|^2dx<\infty.$$ 
\begin{definition}
The quaternionic Segal-Bargmann kernel is defined as $$A(q,x)=\pi^{-\frac{1}{4}}e^{-\frac{1}{2}(q^2+x^2)+\sqrt{2}qx},$$
for every  $q\in \mathbb{H}$ and $x\in \R.$  
\end{definition}
\begin{remark}
The quaternionic Segal-Bargmann kernel $A(q,x)$ is the slice hyperholomorphic extension of the holomorphic Segal-Bargmann kernel $A(z,x)$. Furthermore, it holds that 

$$\displaystyle   A(q,x)=A_q(x) := \sum_{k=0}^{\infty}\psi_k(x)\frac{q^k}{\sqrt{k!}}=\pi^{-1/4} e^{-\frac{1}{2}(x^2+q^2)+\sqrt{2}qx},$$
for every $q\in \mathbb{H}$ and $x\in\mathbb{R}$.
\end{remark}

\begin{definition}
For a given $\varphi\in L^2_\mathbb{H}(\R)$, the quaternionic Segal-Bargmann transform is defined by the following integral expression $$(\mathcal{B}_\mathbb{H}\varphi)(q):=\displaystyle \int_\R A(q,x)\varphi(x)dx.$$ 
\end{definition}

An important result that was proven in \cite{DG2017} is the following:

\begin{theorem}
The quaternionic Segal-Bargmann transform $\mathcal{B}_\mathbb{H}$ defines a surjective isometry mapping the Hilbert space $L^2_\mathbb{H}(\R)$ onto the slice hyperholomorphic Fock space $\mathcal{F}_{Slice}(\mathbb{H})$. Moreover,  $\mathcal{B}_\mathbb{H}$ maps the normalized Hermite functions onto quaternionic monomials. Specifically, we have

\begin{equation}
(\mathcal{B}_{\mathbb{H}}\psi_k)(q)=\frac{q^k}{\sqrt{k!}}, \quad k\geq 0.
\end{equation}
\end{theorem}
\begin{remark}
The quaternionic Fock space $\mathcal{F}_{Slice}(\mathbb{H})$ can be realized as the image of the Hilbert space $L^2_\mathbb{H}(\mathbb{R})$ via the quaternionic Segal-Bargmann transform. Namely, we have 
$$\displaystyle \mathcal{F}_{Slice}(\mathbb{H})=\mathcal{B}_\mathbb{H} \left( L^2_\mathbb{H}(\mathbb{R}) \right)=\left\lbrace{\mathcal{B}_\mathbb{H} \psi \text{ }| \text{ } \psi \in L^2_\mathbb{H}(\mathbb{R}) }\right\rbrace.$$
\end{remark}


Additional results involving the quaternionic Segal-Bargmann transform are proved in \cite{DMD}. For example, by considering the Schwartz space of quaternionic valued functions defined as
$$\displaystyle \mathcal{S}_\mathbb{H}(\mathbb{R}):=\left \lbrace{\psi:\mathbb{R}\longrightarrow\mathbb{H}  \mid \sup_{x\in\mathbb{R}}\left|x^\alpha\frac{d^\beta}{dx^\beta}\psi(x)\right|<\infty},\textbf{ } \forall \alpha,\beta\in \mathbb{N}_0 \right\rbrace,$$

it is possible to study the range of $\mathcal{S}_\mathbb{H}(\mathbb{R})$ under the quaternionic Segal-Bargmann transform $\mathcal{B}_{\mathbb{H}}$, denoted by $\mathcal{SF}(\mathbb{H})$. Indeed, a characterization of $\mathcal{SF}(\mathbb{H})$ is given by the following result:
\begin{theorem}
A function $f(q)=\displaystyle\sum_{k=0}^\infty q^kc_k$ belongs to $\mathcal{SF}(\mathbb{H})$ if and only if $$\sup_{k\in\mathbb{N}_0}|c_k|k^p\sqrt{k!}<\infty, \quad \forall p >0.$$
That is,
$$\displaystyle \mathcal{SF}(\mathbb{H})=\left\lbrace{\displaystyle\sum_{k=0}^\infty q^k c_k \mid  c_k\in\mathbb{H} \text{ and  }\sup_{k\in\mathbb{N}_0}|c_k|k^p\sqrt{k!}<\infty, \quad \forall p >0}\right\rbrace.$$
\end{theorem}

Consider now the position and momentum operators defined by $$X:\varphi\mapsto (X\varphi)(x)=x\varphi(x) \text{ and } P:\varphi\mapsto (P\varphi)(x)=\frac{d}{dx}\varphi(x).$$

Under suitable conditions, it turns out that the quaternionic Segal-Bargmann transform and the slice derivative are related by the following properties (see \cite{DMD}).
\begin{proposition}
\begin{itemize}
\item[(i)] The slice derivative of the quaternionic Segal-Bargmann kernel $A(q,x)$ is given by
$$\partial_S A(q,x)=(-q+\sqrt{2}x) A(q,x), \quad \forall q\in \mathbb H,  \forall x\in \mathbb R.$$

\item[(ii)] Let $\psi \in L^2_{\Hq}(\R)$ such that $\displaystyle X\psi=x\psi$ and $\displaystyle P\psi=\frac{d}{dx}\psi$ belong to  $L^2_{\Hq}(\R)$. Then, the following identities hold:

 $$\displaystyle \left(\partial_S+q\right)\mathcal{B}_{\Hq} \psi=\sqrt{2}\mathcal{B}_{\Hq}X\psi, \text{ and }\mathcal{B}_{\Hq}\left(\frac{X-P}{\sqrt{2}}\right)=q\mathcal{B}_{\Hq}\psi.$$
\end{itemize}
\end{proposition}
\begin{remark}
Another important fact about the Segal-Bargmann transform in the complex case is its connection with the Short-Time Fourier Transform (STFT) with a Gaussian window, as noted in \cite[Proposition 3.4.1]{GK2001}. Based on this observation, the quaternionic Segal-Bargmann transform was used in \cite{DMD} to introduce a quaternionic counterpart of the STFT. Additionally, it was employed to to study a quaternionic analogue of the Gabor space and  its connection with the quaternionic Fock space $\mathcal{F}_{Slice}(\mathbb{H})$.
\end{remark}

\subsection{Gaussian Radial Basis Function (RBF) Kernel of a Quaternionic Variable}\label{Fslicealpha}
The quaternionic slice hyperholomorphic counterpart of the Gaussian RBF kernel has been studied in \cite{DMD3}. To introduce this notion, we briefly recall the definition of the slice hyperholomorphic Fock space, which depends on a parameter $\nu>0$.
\begin{definition}\label{Fslicealpha}
	Let $\nu>0$ be a real parameter. For a given $I \in \mathbb{S}$, the slice hyperholomorphic Fock space is defined as
	$$ \mathcal{F}_{Slice}^{\nu}(\mathbb{H}):= \left \{ f \in \mathcal{SH}(\mathbb{H}) \mid \frac{\nu}{\pi} \int_{\mathbb{C}_I} |f_I(q)| e^{-\nu|q|^2} d \lambda_I(q)<\infty \right\},$$
where $f_I=f_{|\mathbb{C}_I}$ is the restriction of $f$ to $\mathbb{C}_I$, and $ d\lambda_I(q)=dxdy$ is the Lebesgue measure with respect to the variable $q=x+Iy$.	
\end{definition}
\begin{remark}
Note that the case $\nu=1$ corresponds to Definition \ref{SliceFock}. Similarly, the space $ \mathcal{F}_{Slice}^{\nu}(\mathbb{H})$ is a right quaternionic Hilbert space equipped with the inner product
$$ \langle f,g \rangle_{\mathcal{F}_{Slice}^{\nu}(\mathbb{H})}=  \frac{\nu}{\pi}  \int_{\mathbb{C}_I} \overline{g_I(q)} f_I(q) e^{- \nu |q|^2} d \lambda_I(q),$$
for $f$, $g \in \mathcal{F}_{Slice}^{\nu}(\mathbb{H})$. 
\end{remark}

The quaternionic monomials $q^n$ form an orthogonal basis of $ \mathcal{F}_{Slice}^\nu(\mathbb{H})$ with
\begin{equation}
	\label{ort}
	\langle q^m, q^n \rangle_{ \mathcal{F}_{Slice}^\nu(\mathbb{H})}= \frac{m!}{\nu^{m}} \delta_{n,m}.
\end{equation}
The reproducing kernel of $ \mathcal{F}_{Slice}^\nu(\mathbb{H})$ is given by
$$ K_\nu(p,q)=e_{*}^\nu(q\bar{p})= \sum_{n=0}^\infty \frac{\nu^n q^n \bar{p}^n}{n!}.$$

Inspired by the complex case (see \cite{ADCS2022,SC2008,SDC}), we now present the quaternionic Gaussian RBF space in the slice hyperholomorphic setting.
\begin{definition}
Let $\gamma>0$. The slice hyperholomorphic Gaussian RBF space is defined as
$$
\mathcal{H}_{\gamma,I}(\mathbb{H}):= \left\{f \in \mathcal{SH}(\mathbb{H}) \mid  \frac{2}{ \pi \gamma^2} \int_{\mathbb{C}_I} |f_I(q)| e^{\frac{(q- \bar{q})^2}{\gamma^2}} d\lambda_I(q) < \infty \right\}, \quad \forall I \in \mathbb{S}.
$$
\end{definition}
The right $ \mathbb{H}$-vector space $ \mathcal{H}_{\gamma,I}(\mathbb{H})$ is endowed with the inner product
$$ \langle f,g \rangle_{\mathcal{H}_{\gamma,I}(\mathbb{H})}=\frac{2}{ \pi \gamma^2} \int_{\mathbb{C}_I} \overline{g_I(q)} f_I(q)e^{\frac{(q- \bar{q})^2}{\gamma^2}} d \lambda_I(q) ,$$
for $f$, $g \in \mathcal{H}_{\gamma,I}(\mathbb{H})$. 

\begin{remark}
By following similar arguments as in \cite[Theorem 3.1]{DG2017} and using the representation formula, one can show that the slice hyperholomorphic RBF-space does not depend on the choice of the imaginary unit $I \in \mathbb{S}$. Therefore, we can denote the space $\mathcal{H}^I_{\gamma}(\mathbb{H})$ simply by $\mathcal{H}_{\gamma,S}(\mathbb{H})$.
\end{remark}

An orthonormal basis of the slice hyperholomorphic RBF space is given by 
\begin{equation}
\label{basis}
e_{k}^{\gamma}(q)= \sqrt{\frac{2^k}{\gamma^{2k} k!}} q^k e^{-\frac{q^2}{\gamma^2}}, \qquad \gamma >0, k\in \mathbb{N}_0.
\end{equation}

The relation between the spaces $\mathcal{H}_{\gamma,S}(\mathbb{H})$ and $ \mathcal{F}_{Slice}^{\nu}(\mathbb{H})$ when $\nu=\frac{2}{\gamma^2}$ is established via the following result:
\begin{theorem}
\label{one}
Let $\gamma>0$, an entire slice hyperholomorphic function $f: \mathbb{H} \to \mathbb{H}$ belongs to the slice hyperholomorphic RBF space $\mathcal{H}_{\gamma,S}(\mathbb{H})$ if and only if there exists a unique function $g$ in the slice hyperholomorphic Fock space $ \mathcal{F}_{Slice}^{\frac{2}{\gamma^2}}(\mathbb{H})$ such that
$$ f(q)=e^{- \frac{q^2}{\gamma^2}} g(q), \qquad \forall q \in \mathbb{H}.$$
\end{theorem}

\begin{remark}
An isometric isomorphism between the spaces $\mathcal{H}_{\gamma,S}(\mathbb{H})$ and $ \mathcal{F}_{Slice}^{\frac{2}{\gamma^2}}(\mathbb{H})$ is given by
$$ \mathcal{M}^{\gamma^2}[f](q)= e^{\frac{q^2}{\gamma^2}} f(q), \quad f \in \mathcal{H}_{\gamma,S}(\mathbb{H}), \quad q \in \mathbb{H}.$$
\end{remark}
\begin{definition}
\label{ker}
Let $ \gamma >0$. The function
$$ K_{\gamma,S}(q,p)=K^{p}_{\gamma,S}(q):= e^{- \frac{q^2}{\gamma^2}} e_{*}^{\frac{2}{\gamma^2}}(q \bar{p})e^{- \frac{\bar{p}^2}{\gamma^2}},$$
is called the quaternionic slice hyperholomorphic Gaussian RBF kernel.
\end{definition}
\begin{remark}
	The quaternionic RBF kernel $K_{\gamma,S}(q,p)$ is slice hyperholomorphic in the variable $q$ and anti-slice hyperholomorphic in the variable $p$. 
\end{remark}	
\begin{theorem}
The slice hyperholomorphic RBF space $ \mathcal{H}_{\gamma,S}(\mathbb{H})$ is a quaternionic reproducing kernel Hilbert space whose kernel is given by the quaternionic Gaussian RBF kernel $K_{\gamma, S}$.
Moreover, the reproducing kernel property is given by the following integral representation
$$
f(p)=  \frac{2}{ \pi \gamma^2}\int_{\mathbb{C}_I} \overline{K_{\gamma,S}(q,p)} f_I(q) e^{\frac{(q- \bar{q})^2}{\gamma^2}} d \lambda_I(q), \quad q \in \mathbb{H}, \qquad \forall f \in \mathcal{H}_{\gamma,S}(\mathbb{H}).
$$
\end{theorem}

The following properties are satisfied by the quaternionic Gaussian RBF kernel $K_{\gamma,S}$:
\begin{proposition} It holds that 
\begin{enumerate}
\item $\displaystyle K_{\gamma,S}(q,p)= \sum_{k=0}^\infty e_{k}^{\gamma}(q) \overline{e_k^\gamma(p)},$ for every $q,p\in \mathbb{H}$,
\item $ \langle K_{\gamma,S}^q, K_{\gamma,S}^p \rangle_{\mathcal{H}_{\gamma,S}(\mathbb{H})}=K_{\gamma,S}(p,q),$ for every $q,p\in \mathbb{H}$,
\item $ | f(q)| \leq e^{\frac{2}{\gamma^2} y^2} \| f \|_{\mathcal{H}_{\gamma,S}(\mathbb{H})},$ for every  $ f \in \mathcal{H}_{\gamma,S}(\mathbb{H})$ and $q=x+Iy \in \mathbb{C}_I$.
\end{enumerate}
\end{proposition}
\begin{remark}
The quaternionic Segal-Bargmann transform was used in \cite{DMD} to introduce the quaternionic RBF Bargmann transform, which maps $L^2_{\mathbb H}(\mathbb R)$ onto the RBF space $\mathcal{H}_{\gamma,S}(\mathbb{H})$.
\end{remark}

 \subsection{Quaternionic Slice Hyperholomorphic Cholewinski-Fock Space}
The generalized Fock space introduced by Cholewinski in \cite{Cholewinski1984} (see also \cite{SifiSoltani2002}) has been studied in \cite{D2019} in the setting of slice hyperholomorphic functions. To present this extension in this survey, we begin by recalling some basic properties on modified Bessel functions (see \cite{Erdelyi1953,Lebedev1972}).

The modified Bessel function of the first kind is given by $$\mathsf I_\nu(x):=\displaystyle\left(\frac{x}{2}\right)^\nu\sum_{k=0}^\infty \frac{1}{\Gamma(k+1)\Gamma(\nu+k+1)}\left(\frac{x}{2}\right)^{2k}.$$ The Macdonald function is defined as: $$\displaystyle \mathsf K_\nu(x)=\frac{\pi}{2}\frac{\mathsf I_{-\nu}(x)-\mathsf I_{\nu}(x)}{\sin(\nu\pi)}\quad\text{if} \quad \nu \notin \Z,$$
and $$\mathsf K_n(x)=\underset{\nu\rightarrow n}\lim \mathsf K_\nu(x) \quad \text{if} \quad\nu=n\in\Z.$$
The Macdonald function plays an important role in the construction of the quaternionic Cholewinski-Fock space. We recall some important properties of this function that will be useful for this purpose, see \cite{Erdelyi1953,Lebedev1972}. \\
\begin{proposition} \label{Mac}
Let $x>0$ and $\delta,\nu\in \R$ such that $\delta+\nu >0$ and $\delta-\nu >0$. The following identities hold:
\begin{enumerate}
\item $\displaystyle \mathsf K_\nu(x)=\int_0^\infty \exp(-x\cosh t)\cosh(\nu t)dt,$
\item $\displaystyle \mathsf K_{\frac{1}{2}}(x)=\mathsf K_{-\frac{1}{2}}(x)=\sqrt{\frac{\pi}{2x}}e^{-x},$
\item $\displaystyle\int_0^\infty t^{\delta-1}\mathsf K_\nu(t)dt=2^{\delta-2}\Gamma\left(\frac{\delta}{2}+\frac{\nu}{2}\right)\Gamma\left(\frac{\delta}{2}-\frac{\nu}{2}\right).$

\end{enumerate}
\end{proposition}

First, we note that any quaternionic entire slice hyperholomorphic function can be written as $$f=f^e+f^o,$$  where $f^e$ and $f^o$ are the even and odd functions, respectively, given by  $$f^e(q):=\displaystyle\frac{f(q)+f(-q)}{2} \quad \text{and} \quad f^o(q):=\frac{f(q)-f(-q)}{2}.$$
The series expansion theorem of slice hyperholomorphic functions allows us to write $$f(q)=\displaystyle\sum^\infty_{n=0} q^na_n \quad \text{with} \quad a_n\in\Hq,$$ so that $$f^e(q)=\displaystyle\sum_{n=0}^\infty q^{2n}a_{2n} \quad \text{and} \quad f^o(q)=\displaystyle\sum_{n=0}^\infty q^{2n+1}a_{2n+1}.$$
Let $\displaystyle \alpha\geq-\frac{1}{2}$ and let $I$ be any imaginary unit in the sphere $\mathbb{S}$. Then, for $p=x+yI$ in the slice $\C_I$, we consider the probability measure $$\displaystyle d\lambda_{\alpha,I}(p):=\frac{\vert{p}\vert^{2\alpha+2}}{\pi 2^{\alpha}\Gamma(\alpha+1)}\mathsf K_\alpha(\vert{p}\vert^2)d\lambda_I(p),$$
where $\mathsf{K}_\alpha$ is the Macdonald function and $d\lambda_I(p)$ is the usual Lebesgue measure on the slice $\C_I$. Inspired by the classical complex case studied in \cite{SifiSoltani2002}, we introduce the following definition.
\begin{definition}
An entire slice hyperholomorphic function $f:\Hq\longrightarrow \Hq$ belongs to the slice Cholewinski-Fock space (or the generalized slice Fock space) if, for some $I\in\Sq$ the following condition holds:
$$\displaystyle \int_{\C_I}\vert{f^e_I(p)}\vert^2 d \lambda_{\alpha,I}(p)+2(\alpha+1)\int_{\C_I}\vert{f^o_I(p)}\vert^2\vert{p}\vert^{-2}d\lambda_{\alpha+1,I}(p)<\infty.$$
The space containing all such functions will be denoted by $\mathcal{F}_{CS}^\alpha(\Hq)$.
\end{definition}

The quaternionic Cholewinski-Fock space $\CFH$ is equipped with the inner product

$$\displaystyle\scal{f,g}_{\CFH}:=\int_{\C_I}\overline{g^e_I(p)}f^e_I(p)d\lambda_{\alpha,I}(p)+2(\alpha+1)\int_{\C_I}\overline{g^o_I(p)}f^o_I(p)\vert{p}\vert^{-2}d\lambda_{\alpha+1,I}(p),$$

for all $f,g \in \CFH$. \\ By using the representation formula for slice hyperholomorphic functions, it is possible to prove that the definition of the quaternionic space $\CFH$ does not depend on the choice of the imaginary unit $I$. Moreover, one can prove the following results:
\begin{theorem}\label{scal}
Let $f(q)=\displaystyle\sum_{n=0}^\infty q^na_n$ and $g(q)=\displaystyle\sum_{n=0}^\infty q^nb_n$ be two slice hyperholomorphic functions belonging to $\CFH$. Then, we have $$\scal{f,g}_{\CFH}=\displaystyle \sum_{n=0}^\infty \overline{b_n}a_n\beta_n(\alpha),$$
where $$\displaystyle\beta_n(\alpha):=2^n\left\lfloor \frac{n}{2}\right\rfloor !\frac{\displaystyle\Gamma\left(\left\lfloor\frac{n+1}{2} \right\rfloor +\alpha+1\right)}{\Gamma(\alpha+1)}.$$
Here the symbol $\lfloor . \rfloor$ denotes the integer part.
\end{theorem}
For any $n\in\N$, consider the functions $$\phi_n^\alpha(q)=\displaystyle\frac{q^n}{\sqrt{\beta_n(\alpha)}}.$$
\begin{theorem}
The family of functions $(\phi_n^\alpha)_{n \geq 0}$ forms an orthonormal basis of $\CFH$. Moreover, the quaternionc Cholewinski-Fock space is a reproducing kernel Hilbert space with a reproducing kernel is given by the kernel function
$$\displaystyle \mathcal C_\alpha(p,q)=\sum_{n=0}^\infty\frac{p^n\overline{q}^n}{\beta_n(\alpha)},\quad p, q\in \mathbb H, $$
where $$\displaystyle\beta_n(\alpha):=2^n\left\lfloor \frac{n}{2}\right\rfloor !\frac{\displaystyle\Gamma\left(\left\lfloor\frac{n+1}{2} \right\rfloor +\alpha+1\right)}{\Gamma(\alpha+1)}.$$
\end{theorem}

\begin{remark}
Observe that for $\alpha=-\frac{1}{2},$ one can prove that $\CFH$ coincides with the slice hyperholomorphic Fock space presented in Definition \ref{SliceFock}. This holds because $$\displaystyle d\lambda_{-\frac{1}{2},I}(p):=\frac{1}{2\pi}e^{-\vert p \vert^{2}}d\lambda_I(p).$$ Moreover, in this case, the kernel function $\mathcal C_{-\frac{1}{2}}(.,.)$ is precisely the reproducing kernel of $\mathcal{F}_{Slice}(\Hq)$.
\end{remark}
\begin{remark}
A unitary integral transform on the quaternionic Cholewinski-Fock space, along with related operators, was introduced and studied in \cite{D2019} using the generalized Hermite polynomials considered in \cite{Rosenblum1994}.
\end{remark}
\subsection{Quaternionic $p$-Fock spaces and approximation}
The $p$-Fock spaces of slice hyperholomorphic functions were first introduced in \cite{DGS2019} to study approximation results in the quaternionic setting. Similarly to the case of quaternionic Bergman spaces considered in \cite{GS2018} (see also \cite{GS2019}), one can define quaternionic $p$-Fock spaces of the first and second kind to establish quaternionic approximation results in these spaces.
\begin{definition}[The First Kind]\label{First}
Let $0<p<+\infty$ and $0<\alpha<+\infty$. The Fock space of the first kind $\mathcal{F}^{p}_{\alpha}(\mathbb{H})$ is defined as the space of entire slice hyperholomorphic functions $f\in \mathcal{SH}(\mathbb{H})$ such that
$$\|f\|_{p, \alpha}^p:=\displaystyle \frac{\alpha p}{2\pi}\int_{\mathbb{H}}\left|f(q)e^{-\alpha|q|^{2}/2}\right|^{p}d\lambda(q)<+\infty,$$
with $d\lambda(q)$ denotes the Lebesgue volume element in $\mathbb{R}^{4}$.
\end{definition}
\begin{remark}
 Using the same techniques as in the complex case, one may verify that for $1\le p <+\infty$, $\|\cdot\|_{p, \alpha}$ satisfies the properties of a norm.
\end{remark}
To introduce the Fock spaces of the second kind, we need the following definition:
\begin{definition}
For $I\in \mathbb{S}$, $0<\alpha <+\infty$, and $0<p<+\infty$, let us denote
$$
 \|f\|_{p, \alpha, I}^p=\displaystyle\frac{\alpha p}{2\pi}\int_{\mathbb{C}_{I}}\left|f(q)e^{-\alpha|q|^{2}/2}\right|^{p}d \lambda_{I}(q),
$$
 where $d \lambda_{I}(q)$ represents the area measure on $\mathbb{C}_{I}$. The space of all entire functions $f$ satisfying $\|f\|_{p, \alpha, I}<+\infty$ will be denoted by $\mathcal{F}^{p}_{\alpha, I}(\mathbb{H})$.
\end{definition}

We are now in position to introduce the following:
\begin{definition}[The Second Kind]\label{Second}
The Fock space of the second kind $\mathcal{F}^{\alpha, p}_{Slice}(\mathbb{H})$ is defined as the space of all $f\in \mathcal{SH}(\mathbb{H})$ such that for some $I\in \mathbb{S}$, we have $f\in \mathcal{F}^{p}_{\alpha, I}(\mathbb{H})$. To ensure that the norm is independent of the choice of the imaginary unit, we set $$\|f\|_{\mathcal{F}^{\alpha, p}_{Slice}(\mathbb{H})}=\underset{I\in\mathbb{S}} \sup \|f\|_{p, \alpha, I}.$$
\end{definition}
We now state the polynomial approximation result for the quaternionic Fock spaces of the first kind.
\begin{theorem}[Approximation of the First Kind]\label{thm_first_kind} Let $\alpha > 0$ and $0<p<\infty$. The set of all quaternionic polynomials is contained in $\mathcal{F}^{p}_{\alpha}(\mathbb{H})$, and for every $f\in \mathcal{F}^{p}_{\alpha}(\mathbb{H})$, there exists a sequence of quaternionic polynomials $(p_{n})_{n\in \mathbb{N}}$ such that $$\|p_{n}-f\|_{p, \alpha}\to 0 \text{ as } n\to +\infty.$$
\end{theorem}

Next, we present the approximation result for the Fock spaces of the second kind.
\begin{theorem}[Approximation of the Second Kind]\label{thm_second_kind} Let $0<p<+\infty$, $0<\alpha<+\infty$, and $f\in \mathcal{F}^{\alpha, p}_{Slice}(\mathbb{H})$. There exists a sequence of polynomials $(P_{n})_{n\in \mathbb{N}}$ such that for any $I\in \mathbb{S}$, 

$$\|P_{n}-f\|_{p, \alpha, I}\to 0 \text{ as }n\to +\infty.$$
\end{theorem}
\begin{remark}
Note that when $p=2$, the notion of the Fock space of the second kind coincides with the slice hyperholomorphic Fock space with a parameter $\alpha>0$, as presented in Definition \ref{Fslicealpha}. In this case, we obtain a reproducing kernel Hilbert space, with the kernel given by 
$$K_\alpha(r,q):= e_*(\alpha r \overline{q})=\displaystyle \sum_{k=0}^\infty \frac{\alpha^k r^k\overline{q}^k}{k!}, \quad r, q\in \mathbb H.$$
\end{remark}

Let $\mathcal{R}$ denote the set of functions of the form
$$
f(r) = \sum_{k=1}^n K_\alpha(r, a_k) \, b_k, \quad \text{for all } r \in \mathbb{H},
$$
where $n \in \mathbb{N}$ and $(a_k)_{k=1,\cdots, n}$ and $(b_k)_{k=1,\cdots, n}$ are quaternions. 

In other words,
$$
\mathcal{R} := \left\{ f = \sum_{k=1}^n K_\alpha(\cdot, a_k) \, b_k \ \bigg|\  n \in \mathbb{N},\ a_k, b_k \in \mathbb{H} \right\}.
$$

An interesting consequence of Theorem \ref{thm_second_kind} in the second-kind theory leads to the following result:
\begin{theorem}\label{RK2} Let $\alpha > 0$ and $0<p<\infty$. The set $\mathcal{R}$ is dense in the quaternionic Fock spaces of the second kind $\mathcal{F}^{\alpha, p}_{Slice}(\mathbb{H})$.
\end{theorem}

\begin{remark}\label{ObsAppro}
New approximation results for complex and quaternionic polyanalytic functions have been studied in \cite{GS2022, GS2023, GS2024}. In particular, the density of polyanalytic polynomials in complex and quaternionic polyanalytic weighted Bergman spaces was considered in \cite{GS2023}, while density results in quaternionic polyanalytic Fock spaces were addressed in \cite{GS2024}.
\end{remark}

\section{Fock spaces of slice polyanalytic functions}
\setcounter{equation}{0}
The Fock space of slice polyanalytic functions of a quaternionic variable can be introduced as follows:
\begin{definition}
Let $n\in \mathbb N$ and $I\in\Sq$ we define the slice polyanalytic Fock space as
$$\displaystyle \mathcal{F}^{n}_{I}(\Hq):=\left\lbrace{f\in\mathcal{SP}_n(\Hq) \mid  ||f||^{2}_{\mathcal{F}^{n}_{I}(\Hq)}:= \int_{\C_I}\vert{f_I(p)}\vert^2e^{-\vert{p}\vert^2}d\lambda_I(p)<\infty }\right\rbrace.$$
\end{definition}
The space $\mathcal{F}^{n}_{I}(\Hq)$ is a right quaternionic Hilbert space with respect to the following inner product
$$\displaystyle\scal{f,g}_{\mathcal{F}^{n}_{I}(\Hq)}=\int_{\C_I}\overline{g_I(p)}f_I(p)e^{-\vert{p}\vert^2}d\lambda_I(p).$$
We make the following observation:
\begin{proposition}
Let $f\in\mathcal{SP}_n(\Hq)$ and $I,J\in\Sq$ be two imaginary units. Then, we have the following $$ \frac{1}{2}\Vert{f}\Vert_{\mathcal{F}^{n}_{I}(\Hq)}\leq{\Vert{f}\Vert_{\mathcal{F}^{n}_{J}(\Hq)}}\leq 2\Vert{f}\Vert_{\mathcal{F}^{n}_{I}(\Hq)}.
  $$
\end{proposition}
\begin{remark}
The slice polyanalytic Fock space does not depend on the choice of the imaginary unit since the norms are equivalent. In other words, a function $f$ belongs to $\mathcal{F}^{n}_{I}(\Hq)$ if and only if it belongs to $\mathcal{F}^{n}_{J}(\Hq)$.  Therefore, we denote this space by $\mathcal{F}^{n}_{Slice}(\Hq).$
\end{remark}

Let $q\in\Hq$ and consider the evaluation mapping $$\Lambda_q: \mathcal{F}^{n}_{Slice}(\Hq)\longrightarrow \Hq,  f\mapsto\Lambda_q(f)=f(q).$$ Then, we have the following estimate.
\begin{proposition} \label{Growth}
Let $f\in \mathcal{F}^{n}_{Slice}(\Hq)$ and $q\in\Hq$. Then, $$|\Lambda_q(f)|\leq \sqrt{n} e^{\frac{|q|^2}{2}}||f||.$$
\end{proposition}

Proposition \ref{Growth} shows that all the evaluation mappings on $\mathcal{F}^{n}_{Slice}(\Hq)$ are continuous. Then, the Riesz representation theorem for quaternionic right-linear Hilbert spaces, see \cite{BDS1982} implies that for any $q\in\Hq$ there exists a unique function $K_N^q\in \mathcal{F}^{n}_{Slice}(\Hq)$ such that for any $f\in \mathcal{F}^{n}_{Slice}(\Hq)$ we have $$f(q)=\scal{f,K^q_n}_{\mathcal{F}^{n}_{Slice}(\Hq)}.$$
Let $J\in\Sq$ and $r\in\C_J$, then for $q=x+Iy$ and $z=x+Jy$ the corresponding reproducing kernel is obtained by extending the complex polyanalytic Fock kernel to the slice polyanalytic setting. It is given by the following $$K_n:\Hq\times\Hq\longrightarrow \Hq$$
$$\displaystyle K_n(q,r):=\frac{1}{2}\left[K_n(z,r)+K_n(\overline{z},r)\right]+I\frac{J}{2}\left[K_n(\overline{z},r)-K_n(z,r)\right].$$

Let us consider the first non-trivial example.
\begin{example}
For $n=2$, we have
$$\displaystyle K_2(q,r)=ext[e^{z\overline{r}}(2-\vert{z-r}\vert^2)](q)$$
Then, using the quaternionic $*$-product of slice functions (see \cite{ColomboSabadiniStruppa2016}), it follows that $$K_2(q,r)=e_*(q\overline{r})*(2-(q\bar{q}-q\bar{r}-\bar{q}r+r\bar{r})).$$
\end{example}

Inspired by the previous example, we can state the following general result:
\begin{theorem}\label{KerFock}
The set $\mathcal{F}^{n}_{Slice}(\Hq)$ is a right quaternionic reproducing kernel Hilbert space whose reproducing kernel function is given by the following formula:

$$\displaystyle K_n(q,r)=e_*(q\overline{r})*\varphi_n(q,r),\quad  q,r\in\Hq,$$
 where the $*$-multiplication is in the variable $q$ and 
 
 $$\varphi_n(q,r)=\sum_{k=0}^{n-1} (-1)^k {n\choose{k+1}}\frac{1}{k!}(\bar q q-q \bar r -\bar q r + \bar r r)^{k*}.$$
\end{theorem}
\begin{remark}
The reproducing kernel of the quaternionic slice polyanalytic Fock space $\mathcal{F}^{n}_{Slice}(\Hq)$ can be expressed in terms of the generalized Laguerre polynomials. Indeed, if we consider the $*$-power in \eqref{Laguerre} and define
$$x^{k*}:=(\bar q q-q \bar r -\bar q r + \bar r r)^{k*},$$
then we have 
 $$\displaystyle K_n(q,r)=e_*(q\overline{r})* L^{1}_{n-1}(\bar q q-q \bar r -\bar q r + \bar r r), \quad q, r\in \mathbb{H}.$$
 \end{remark}
\begin{remark}
For $n=1$, the space $\mathcal{F}^n_{Slice}(\Hq)$ is exactly the slice hyperholomorphic Fock space, and the reproducing kernel obtained in Theorem \ref{KerFock} corresponds to the kernel function given in \eqref{*Fockkernel}.
\end{remark}

The quaternionic true slice polyanalytic Fock space and the associated true polyanalytic Bargmann transform were considered in \cite{DMD2} to study the Short-Time Fourier Transform(STFT) and related topics in the slice polyanalytic context. We refer to \cite{Asampling, AbreuFeichtinger2014, Vas} for the study of the true polyanalytic Fock spaces in the classical case of one complex variable. We briefly recall the definition of the quaternionic true polyanalytic Fock space, given by
\begin{definition}
\label{FT}
A function $ f: \mathbb{H} \to \mathbb{H}$ belongs to the quaternionic true polyanalytic Fock space $ \mathcal{F}_T^n (\mathbb{H})$ if and only if
$$\displaystyle \int_{\mathbb{C}_I}|f_I(q)|^2 e^{-|q|^2} \, d \lambda_I(q) < \infty,$$
and there exists an entire slice hyperholomorphic function $g$ such that
$$f(q)= c(n) e^{|q|^2} \partial_s^{n-1}(e^{-|q|^2} g(q)),$$
 for every $q\in \mathbb{H}$, where $\partial_s $ is the slice derivative and $c(n)$ is a constant depending on the order of polyanalyticity $n$.

\end{definition}

The relation between the quaternionic slice polyanalytic Fock space and the quaternionic true polyanalytic Fock space is given by the following result:
\begin{theorem}
\label{sum}
The quaternionic polyanalytic Fock space $ \mathcal{F}_{Slice}^{n}(\mathbb{H})$ is the direct sum of true polynalytic Fock spaces $\mathcal{F}_T^{\ell}(\mathbb{H})$, $\ell=0,...,n-1$, i.e.,
\label{R3}
$$ \mathcal{F}_{Slice}^{n}(\mathbb{H})= \bigoplus_{\ell=0}^{n-1} \mathcal{F}_T^{\ell}(\mathbb{H}).$$
\end{theorem}

Finally, we present the quaterionic slice polyanalytic counterparts of the $p$-Fock spaces. As discussed in Remark \ref{ObsAppro}, new approximation results have been obtained concerning function spaces of slice polyanalytic functions in the quaternionic variable. For instance, inspired by Definition \ref{polyFp} of the complex polyanalytic $p$-Fock spaces, the authors of \cite{GS2024} introduced the quaternionic counterparts of slice polyanalytic $p$-Fock spaces of the first and second kinds to study quaternionic approximation on these spaces. Specifically, the following definitions can be considered:

\begin{definition}[Poly-Fock First Kind]
Let $0<p<+\infty, n\in \mathbb{N},$ and $0<\alpha<+\infty$. The polyanalytic quaternionic $p$-Fock space of the first kind, denoted by $\mathcal{F}^{p}_{\alpha,n}(\mathbb{H})$, is defined as the space of all slice polyanalytic functions of order $n$ on $\mathbb H$ such that
$$\|f\|_{p, \alpha}^p:=\displaystyle \frac{\alpha p}{2\pi}\int_{\mathbb{H}}\left|f(q)e^{-\alpha|q|^{2}/2}\right|^{p}d\lambda(q)<+\infty,$$
with $d\lambda(q)$ denotes the Lebesgue volume element in $\mathbb{R}^{4}$.
\end{definition}

\begin{definition}[Poly-Fock Second Kind]
Let $0<p<+\infty, n\in \mathbb{N},$ and $0<\alpha<+\infty$. The polyanalytic quaternionic  $p$-Fock space of the second kind, denoted
by $\mathcal{F}^{\alpha, p}_{Slice, n}(\mathbb{H})$, is defined as the space of all slice polyanalytic functions of order $n$ on $\mathbb H$ such that for some $I\in \mathbb{S}$, we have
$$\|f\|_{p, \alpha,I}^p:=\displaystyle \frac{\alpha p}{2\pi}\int_{\mathbb{C}_I}\left|f(q)e^{-\alpha|q|^{2}/2}\right|^{p}d\lambda_I(q)<+\infty.$$
\end{definition}

\begin{remark}
Note that for $n=1$, we recover the slice hyperholomorphic  $p$-Fock spaces of the first and second kinds introduced in Definitions \ref{First} and \ref{Second}. 
\end{remark}

\section{The Fueter mapping theorem and quaternionic Fock spaces}
\setcounter{equation}{0}
In this section, we discuss how the Fueter mapping theorem and its polyanalytic extension connect with the quaternionic slice hyperholomorphic and slice polyanalytic Fock spaces.
\subsection{The Fock-Fueter Space}
\setcounter{equation}{0}
Note that the original version of the Fueter mapping theorem can be found in \cite{Fueter1}. For our purposes, we recall the statement of the Fueter mapping theorem as presented in \cite{ColomboSabadiniSommen2010}:
\begin{theorem}\label{FMthm}
Let $f:U\subset \Hq\longrightarrow \Hq$ be a slice hyperholomorphic function of the form $$f(x+yI)=\alpha(x,y)+I\beta(x,y),$$where $\alpha(x,y)$ and $\beta(x,y)$ are quaternionic-valued functions satisfying the compatibility conditions $$\alpha(x,-y)=\alpha(x,y), \quad \beta(x,-y)=-\beta(x,y),$$ and the Cauchy-Riemann system. Then $$\overset{\sim}f(x_0+\vec{q}\,)=\displaystyle \Delta\left(\alpha(x_0,\vert{\vec{q}\,}\vert)+\frac{\vec{q}\,}{{|\vec{q}\,|}}\beta(x_0,\vert{\vec{q}\,}\vert)\right)$$
defines a Fueter regular function.
\end{theorem}

\begin{remark}
We denote the Fueter mapping by $$\tau:\mathcal{SR}(U)\rightarrow\mathcal{FR}(U), \text{ } f\longmapsto \tau(f)=\Delta f.$$
\end{remark}

The slice hyperholomorphic functions can be represented as a series expansion of quaternionic monomials as stated in Theorem \ref{expansion}. Thus, a natural question is how the Fueter mapping theorem can be applied to the quaternionic monomials, which form a building block four understanding slice hyperholomorphic functions. Basically, the authors of \cite{DKS2019} computed the action of the Fueter map $\tau$ on the quaternionic monomials $q^k$, leading to the following result:

\begin{lemma}[\cite{DKS2019}] \label{FueterAct}
The action of the Fueter map $\tau$ on the quaternionic monomials $\lbrace q^k \mid k\in \mathbb N_0\rbrace $ is given by
$$\tau(1)=\tau(q)=0,\quad \tau(q^k):=\Delta (q^k)=-4\displaystyle\sum_{s=1}^{k-1}(k-s)q^{k-s-1}\overline{q}^{s-1}, \quad k\geq 2.$$
\end{lemma}

Let $k\geq 0$ and define the following family of Fueter regular functions \begin{equation} \label{Qk1}
\displaystyle Q_k(q):=-\frac{\tau(q^{k+2})}{2(k+1)(k+2)}=-\frac{\Delta(q^{k+2})}{2(k+1)(k+2)},\quad k\in \mathbb N_0.
\end{equation} 

By comparing the result obtained in Theorem \ref{FueterAct} with the coefficients used in formulas 5 and 6 of \cite{CMF2011} (see also \cite{MF2007}), we make the following observation:

\begin{theorem}

 The Fueter regular functions $(Q_k)_{k\in \mathbb N_0}$ can be written as

 \begin{equation} \label{Qk2}
\displaystyle Q_k(q)=Q_k(q,\overline{q})=\sum_{j=0}^k T^k_jq^{k-j}\overline{q}^j, \quad q\in \mathbb H,  k \in \mathbb N_0,
\end{equation}

where $$\displaystyle T^k_j:=T^{k}_{j}(3)=\frac{k!}{(3)_k}\frac{(2)_{k-j}(1)_j}{(k-j)!j!}=\frac{2(k-j+1)}{(k+1)(k+2)},$$ and $(a)_n=a(a+1)...(a+n-1)$ is the Pochhammer symbol. Moreover, the Fueter regular polynomials $Q_k$ satisfy the Appell property with respect to the hypercomplex derivative, i.e., 
\begin{equation}\label{Appell}
\displaystyle\frac{1}{2}\overline{D}Q_k:=\displaystyle \frac{1}{2}\left(\frac{\partial}{\partial x_0}-i\frac{\partial}{\partial x_1}-j\frac{\partial}{\partial x_2}-k\frac{\partial}{\partial x_3}\right) Q_k =kQ_{k-1}, \quad k\in \mathbb N.
\end{equation}
\end{theorem}
\begin{remark}
We observe that $$\displaystyle \Delta(q^k)=-2(k-1)kQ_{k-2}(q), \quad  k\geq 2.$$
\end{remark}

\begin{definition}
The Fueter regular functions defined as $$\displaystyle Q_k(q):=-\frac{\tau(q^{k+2})}{2(k+1)(k+2)}=\sum_{j=0}^k T^k_jq^{k-j}\overline{q}^j, \quad q\in \mathbb H,  k \in \mathbb N_0,$$

are called Clifford-Appell polynomials.
\end{definition}

\begin{remark}
By applying the Fueter mapping theorem, \cite[Theorem~7.7]{CDS2025} established an extension of the Clifford-Appell polynomials $Q_k(q)$, yielding new forms of Taylor and Laurent series expansions in quaternionic analysis.
\end{remark}

Consider now the following space: 
\begin{definition}\label{Fock-FueterSpace}
The range of the quaternionic slice hyperholomorphic Fock space via the Fueter map $\tau$ is defined by $$ \displaystyle \mathcal{A}(\Hq):=\left\lbrace \tau (f)\mid f\in \mathcal{F}_{Slice}(\Hq)\right\rbrace.$$
The space $\mathcal{A}(\Hq)$ is called the Fock-Fueter space.
\end{definition}
One can prove the following analytic characterization of the space $\mathcal{A}(\Hq)$:

\begin{theorem}[\cite{DKS2019}]  \label{AHcar}
A Fueter regular function $f$ belongs to $ \mathcal{A}(\Hq)$ if and only if
 $$\displaystyle f=\sum_{k=0}^\infty Q_k\alpha_k,$$ 
 
 where $(\alpha_k)_{k\geq 0}$ are quaternionic coefficients, and
$$\displaystyle \sum_{k=0}^{\infty}\frac{k!}{(k+1)(k+2)}|\alpha_k|^2<\infty.$$
 \end{theorem}
 
 \begin{theorem}[\cite{DKS2019}]  
For every $k\in \mathbb N_0,$ consider the quaternionic Fueter regular polynomials defined by \begin{equation}
T_k(q)=\displaystyle\sqrt{\frac{(k+1)(k+2)}{k!}}Q_k(q), \quad q\in\Hq.
\end{equation}
The reproducing kernel of the Fock-Fueter space $\mathcal{A}(\Hq)$ is given by
\begin{equation}\label{circle}
G(p,q)=G_q(p):=\displaystyle \sum_{k=0}^{\infty}T_k(p)\overline{T_{k}(q)}, \quad p,q\in \mathbb H.
\end{equation}
\end{theorem}

A Bargmann-Fock-Fueter transfrom was studied in \cite{DKS2019} by considering the commutative diagram $$\mathcal{S}_\Hq: \xymatrix{
    L^2_\Hq(\R) \ar[r] \ar[d]_{\mathcal{B}_\Hq} & \mathcal{A}(\Hq)  \\ \mathcal{F}_{Slice}(\Hq) \ar[r]_{Id} & \mathcal{SR}(\Hq) \ar[u]_{\tau}
  }$$
  so that $$\mathcal{S}_\Hq:=\tau\circ Id \circ \mathcal{B}_\Hq.$$
   
   \begin{remark}
Various extensions and several results based on the idea of using this Fueter mapping approach to study function spaces over the quaternions have been developed (see \cite{ADCS2021, ADSMN, DMDG1, DMDG2}). 
 The Clifford-Appell polynomials $(Q_k)_{k\in\mathbb N_0}$ also appear in \cite{DP2023} while studying a polyanalytic functional calculus on the $S$-spectrum, as well as in \cite{KP2022} in the context of eigenvalue problems associated with axially monogenic functions.
 \end{remark}

\subsection{The Poly-Fock-Fueter Space}

The Fueter mapping theorem admits two possible extensions in the polyanalytic setting, which were studied in \cite{ADSP2020} using the quaternionic global operator and a slice polyanalytic Cauchy formula. \\ 

Using the observation in Remark \ref{GVR}, we obtain the first polyanalytic extension of the Fueter mapping, stated in the following result:
\begin{theorem}[Poly-Fueter Mapping Theorem I]
Let $\Omega$ be an axially symmetric slice domain in $\Hq$, and let $f:\Omega\longrightarrow \Hq$ a slice polyanalytic function of order $n\geq 1$. Then, the function defined as
$$\displaystyle \tau_n(f)=\Delta \circ V^{n-1}f, $$
where $V$ is the quaternionic modified global operator, belongs to the kernel of the Cauchy-Fueter operator $D$.
\end{theorem}
The second poly-Fueter theorem is based on the poyanalytic decomposition and is stated as follows:
\begin{theorem}[Poly-Fueter Mapping Theorem II]
Let $\Omega\subseteq\mathbb H$ be an axially symmetric slice domain, and let $f:\Omega\longrightarrow \Hq$ be a slice polyanalytic function of order $n\geq 1$. Assume that $f$ admits the poly-decomposition

$$\displaystyle f=\sum_{k=0}^{n-1}\overline{q}^kf_k,$$ where $f_0,...,f_{n-1}\in \mathcal{SH}(\Omega).$
Then, the function defined by
\begin{equation}
\displaystyle \mathcal{C}_n(f)=\sum_{k=0}^{n-1}x_0^k\Delta f_k
\end{equation}

is a poly-Fueter regular function of order $n$.
\end{theorem}

The first poly-Fueter map $\tau_n$ takes the space of slice polyanalytic functions of order $n\in \mathbb N$ into the space of Cauchy-Fueter regular functions $\mathcal{FR}(\Omega)$. In contrast,
the second poly-Fueter map $\mathcal{C}_n$ allows for the construction of poly-Fueter regular functions starting from slice polyanalytic functions of the same order. 

\begin{remark}
The two poly-Fueter mapping theorems are related through the commutative diagram:
 \[
\xymatrix{
\mathcal{SP}_n \ar[r]^{\tau_n} \ar[d]_{\mathcal{C}_n}& \mathcal{FR}\\
\mathcal{FR}_n\ar[ru]_{(2D)^{n-1}} }
\]
such that  $$\displaystyle \tau_n:=(2D)^{n-1}\circ \mathcal{C}_n,$$

\end{remark}

The poly-Fueter mapping theorems were employed in \cite{DMD2022} to extend the Clifford-Appell polynomials $(Q_k)_{k\in \mathbb N_0}$ to the case of poly-Fueter regular functions, leading to the so-called poly-Fock-Fueter space, which extends the Fock-Fueter space $\mathcal{A}(\mathbb H)$ presented in Definition \ref{Fock-FueterSpace}. For instance, by considering the action of the second poly-Fueter map $\mathcal{C}_{n}$ on the building blocks $\overline{q}^kq^j$, one can introduce the following definition:

\begin{definition}
Let $n, s\in\mathbb N_0$. The generalized Appell polyanalytic polynomials are defined as
\begin{equation}
\label{polmal}
\mathcal{M}_{k,s}(q,\overline{q}):=x_0^kQ_s(q,\overline{q}),\quad k=0,...,n.
\end{equation}
\end{definition}

An extension of the classical Appell property \eqref{Appell} is given by the following result: 

\begin{theorem}
\label{appe1}
For any fixed $n\in \mathbb N_0$ and $s \geq n+1$, we have
		
$$\overline{D}^{n+1}\mathcal{M}_{n,s}(q,\overline{q})=\displaystyle\sum_{j=1}^{n+1} 2^j {n+1\choose j}\frac{n!s!}{(j-1)!(s-j)!}\mathcal{M}_{j-1,s-j}(q,\overline{q}). $$
\end{theorem}

\begin{remark}

The polynomials $\mathcal{M}_{k,s}$ are poly-Fueter regular of order $n+1$ for all $k=0,\cdots, n$. Moreover, we observe that Proposition \ref{appe1} extends the classical Appell property satisfied by the Clifford-Appell polynomials $(Q_s)_{s\in \mathbb N_0}$. Indeed, for $n=0$ in Theorem \ref{appe1} and by noting that $\mathcal{M}_{0,s}(q,\bar{q})=Q_s(q,\bar{q})$, we get $$\overline{D}Q_s(q,\bar{q})=2\frac{s!}{(s-1)!}\mathcal{M}_{0,s-1}(q,\overline{q})=2sQ_{s-1}(q,\bar{q}),\quad s\geq 1.$$
\end{remark}

We now introduce the following definition:

\begin{definition}
 The true poly-Fock-Fueter space $\mathcal{A}_{n+1}(\mathbb{H})$ is defined as the range of the quaternionic true slice polyanalytic Fock space $\mathcal{F}_{n+1}^T(\mathbb{H})$ via the poly-Fueter mapping $\mathcal{C}_{n+1}$. Specifically, we have
\begin{equation}
\label{ff}
\displaystyle \mathcal{A}_{n+1}(\mathbb{H}):=\left\lbrace{\mathcal{C}_{n+1}(f)\mid f\in\mathcal{F}_{n+1}^T(\mathbb{H})}\right\rbrace.
\end{equation}
\end{definition}
We state the following characterization of the ture poly-Fock-Fueter space $\mathcal{A}_{n+1}(\mathbb{H})$:
\begin{theorem}\label{Space5}
Let $n\in \mathbb N_0$. A poly-Fueter regular function $f$ belongs to $ \mathcal{A}_{n+1}(\Hq)$ if and only if
 $$\displaystyle f(q,\overline{q})=\sum_{h=0}^{\infty}\sum_{s=0}^n\mathcal{M}_{n-s,h+n-s}(q,\bar{q})\beta_{h,s}, $$ 
 
 where $(\beta_{h,s})_{h \geq 0, \, 0 \leq s \leq n}$ are quaternionic coefficients, and
$$\displaystyle \sum_{h=0}^{\infty}\sum_{s=0}^n\frac{[(h+n-s)!]^2(s!)^2[(n-s)!]^2}{(h+n+2)!}|\beta_{h,s}|^2 < \infty.$$
\end{theorem}
 \begin{remark}
If we take $n=0$ in Theorem \ref{Space5}, we recover the Fock-Fueter space $\mathcal{A}(\mathbb{H})$ considered in Theorem \ref{AHcar}. Indeed, if $n=0$ we observe that
$$\mathcal{A}_1(\mathbb{H})=\left\lbrace \sum_{h=0}^{\infty}\mathcal{M}_{0,h}(q,\bar{q})\beta_{h,0}\mid \quad (\beta_{h,0})_{h\geq 0}\subset \mathbb{H}, \quad \sum_{h=0}^{\infty}\frac{(h!)^2}{(h+2)!}|\beta_{h,0}|^2<\infty \right\rbrace.$$
This observation holds because $\mathcal{M}_{0,h}(q, \bar{q})=Q_h(q, \bar{q})$ and $$ \sum_{h=0}^{\infty}\frac{(h!)^2}{(h+2)(h+1) h!}=\sum_{h=0}^{\infty}\frac{h!}{(h+2)(h+1)}.$$
As a result, we obtain  
$$\mathcal{A}_1(\mathbb{H})=\mathcal{A}(\mathbb{H}).$$
\end{remark}

Finally, we introduce the following:
\begin{definition}
The poly-Fock-Fueter space $ \mathfrak{A}_{N}(\mathbb{H})$ is defined as the direct sum of the true poly-Fock-Fueter spaces $\mathcal{A}_{n+1}(\mathbb{H})$. Specifically, we have 
$$ \mathfrak{A}_{N}(\mathbb{H})= \bigoplus_{n=0}^N \mathcal{A}_{n+1}(\mathbb{H}).$$
\end{definition}
We now state the following characterization for the poly-Fock-Fueter space $ \mathfrak{A}_{N}(\mathbb{H})$:
\begin{theorem}
Let $N \in \mathbb N_0$. Then, the following characterization holds

$$
\begin{aligned}
\mathfrak{A}_{N}(\mathbb{H}) = \Bigg\{ 
& \sum_{n=0}^{N} \sum_{h=0}^{\infty} \sum_{s=0}^{n} \mathcal{M}_{n-s, h+n-s}(q, \bar{q}) \, \beta_{n,h,s} \ \Bigg| \\
& \sum_{n=0}^{N} \sum_{h=0}^{\infty} \sum_{s=0}^{n} 
\frac{[(h+n-s)!]^2 (s!)^2 [(n-s)!]^2}{(n+1)! (h+n+2)!} \left|\beta_{n,h,s}\right|^2 < \infty 
\Bigg\}.
\end{aligned}
$$

\end{theorem}

\section{Conclusion}
\setcounter{equation}{0}

In this chapter, we reviewed fundamental concepts of Fock spaces in one complex variable and discussed several extensions of these spaces in hypercomplex analysis. Firstly, we recalled standard notations of quaternions and presented two prominent function theories of a quaternionic variable. Specifically, we considered the recent theory of slice hyperholomorphic (and polyanalytic) functions and the classical theory of Fueter (and polyanalytic) regular functions. Secondly, we explored various extensions of Fock spaces and related topics studied in the last few years in the context of slice hyperholomorphic and slice polyanalytic functions of a quaternionic variable. Finally, we established links and connections relating these hypercomplex Fock spaces in the slice theory to Fueter regular and poly-Fueter regular functions, based on the well-known Fueter mapping theorem and its polyanalytic extensions. \\ \\
An interseting direction for future research is to investigate a quaternionic counterpart of the Weyl operator in the slice hyperholomorphic (and polyanalytic) setting, and to compare the Fock-Fueter space with the monogenic Fock space introduced by Cnops and Kisil in \cite{CK1999}(see also \cite{GLQ2009}). Another promising research avenue is to futher develop the theory of Fock spaces initiated in \cite{XS2024} using generalized partial-slice monogenic functions. It is also woth noting that Fock spaces and Segal-Bargmann transforms have been studied in hypercomplex analysis from various perspectives, many of which remain to be explored. For instance, connections with coherent states in Clifford analysis are discussed in \cite{DMNQ2018,KMNQ2016,MNQ2017}. A different construction of the slice Fock space, compatible with the slice Fourier transform, is described in \cite{CD2017}. Moreover, several Bargmann-type transfoms in the monogenic setting have been developed in \cite{BESC2019, CSS2017, DKS2019, PSS2017}, and a Clifford short-time Fourier transform was considered in \cite{DMA2022}.
\section*{Acknowledgment} 
This research  is supported by the Research Foundation–Flanders (FWO) under Grant No 1268123N.

\end{document}